\documentclass[final, 3p, times]{elsarticle}
\usepackage{amsmath,amssymb,amsfonts}
\usepackage{graphicx}
\usepackage{textcomp}
\usepackage{caption}
\usepackage{float}
\usepackage{subcaption}
\usepackage{dsfont}
\usepackage{optidef}
 
\usepackage{amsthm}
\usepackage{tikz}
\usetikzlibrary{decorations.pathreplacing}
\usepackage{multicol}
\usepackage{cuted}
\usepackage{multirow}
\usepackage{subcaption}
\newtheorem{theorem}{Theorem}
\newtheorem{lemma}[theorem]{Lemma}
\newtheorem{proposition}[theorem]{Proposition}
\newtheorem{definition}{Definition}
\newtheorem{problem}{Problem}
\newtheorem{assumption}{Assumption}
\newtheorem{remark}{Remark}[section]

\newtheorem{corollary}{Corollary}[theorem]
\usepackage{xcolor}
\usepackage{cite}
\usepackage[colorlinks=true,allcolors=steelblue]{hyperref}
\usepackage{comment}
\usepackage{algorithm}
\usepackage{algpseudocode}
\algnewcommand\algorithmicforeach{\textbf{for each}}
\algdef{S}[FOR]{ForEach}[1]{\algorithmicforeach\ #1\ \algorithmicdo}
\algrenewcommand\algorithmiccomment[1]{\textcolor{blue}{// {\texttt{#1}}}}
\usepackage{scalerel}
\usepackage{bm}
\usepackage{pgfplots}
\pgfplotsset{compat=1.18}
\usepackage{grffile}
\pgfplotsset{compat=newest}
\usetikzlibrary{plotmarks}
\usetikzlibrary{arrows.meta}
\usepgfplotslibrary{patchplots}
\definecolor{steelblue}{RGB}{70,130,180}

\usepackage{lineno}

\journal{Stochastic Processes and their Applications}

\begin{document}
\begin{frontmatter}
\title{Temporal-Logical Stochastic Processes}
\author{Venkatraman Renganathan} 
\affiliation{
organization={Department of Artificial Intelligence, Indian Institute of Technology Hyderabad, Kandi},
state={Telangana},
country={India}}

\begin{abstract}
We develop a measure-theoretic framework for Temporal-Logical Stochastic Processes (TLSP), in which temporal-logic specifications induce time-indexed stochastic observables capturing the robustness of temporal satisfaction. Building upon existing probabilistic interpretations of temporal-logic robustness, we systematically investigate the stochastic-process structure of these observables. We establish their measurability, shifted-filtration adaptedness, and sample-path regularity under suitable assumptions, and characterize the action of temporal operators on stochastic processes. We derive pathwise perturbation bounds and quantify the sensitivity of temporal satisfaction probabilities through the probability mass concentrated near the logical boundary. A topology on the observable-process space and a semantic pseudometric on temporal formulas are introduced to establish convergence of temporal-logical observables under perturbations of formula parameters and temporal intervals. We further characterize conditions ensuring weak convergence of observable-process distributions and convergence of temporal satisfaction probabilities. Finally, we establish a process-level large-deviation principle for temporal-logical observables, providing a variational characterization of rare temporal events through minimum-action trajectories. The resulting framework connects temporal-logic semantics with stochastic-process topology, stability, convergence, and rare-event asymptotics, providing a mathematical foundation for the analysis of temporally constrained stochastic systems.
\end{abstract}


\end{frontmatter}


\section{Introduction}

Stochastic processes provide a fundamental mathematical framework for modeling
time-evolving uncertainty and play a central role in probability theory,
stochastic analysis, and their applications in engineering, finance, and the
natural sciences \citep{oksendal2003stochastic, karatzas1991brownian,
billingsley1995probability, durrett2019probability}. Classical developments
in this area have focused on structural properties such as martingales,
Markov processes, and long-term behavior, forming the backbone of modern
probabilistic analysis \citep{revuz1999continuous, ethier2009markov,
shiryaev1996probability}. In parallel, temporal logic has emerged as a powerful formalism for specifying
time-dependent properties of dynamical systems. Originally introduced in
computer science to reason about program behavior \citep{pnueli1977temporal},
temporal logic has since become a central tool in formal verification,
model checking, and the analysis of cyber-physical systems
\citep{baier2008principles}. It enables concise descriptions of spatio-temporal and event based specifications to be respected by a system.

In recent years, there has been increasing interest in combining stochastic
systems with temporal logic. This is motivated by applications in which
systems evolve under uncertainty while being required to satisfy complex
time-dependent specifications. Examples include autonomous systems operating
in uncertain environments, financial systems subject to risk constraints,
and engineered systems with safety-critical requirements. Several frameworks
have been developed to address such problems, including probabilistic model
checking \citep{kwiatkowska2007stochastic, baier2008principles}, stochastic
reachability analysis \citep{abate2008probabilistic}, and signal temporal
logic for real-valued signals \citep{maler2004monitoring, donze2010robust}. A key development in this direction has been the introduction of quantitative
semantics for temporal logic, particularly in the context of signal temporal
logic. Instead of assigning a binary truth value to a temporal formula, these
approaches associate a real-valued robustness measure that quantifies the
degree to which a trajectory satisfies or violates a specification
\citep{donze2010robust, fasoni2020robustness}. This has enabled the use of
optimization and statistical techniques for analysing temporal behavior of stochastic systems.

Despite the diversity of approaches, the existing literature shares a common structural paradigm: temporal logic is treated as an \emph{external specification} imposed on a stochastic system. In existing paradigm, a stochastic process is first defined, typically as a probability measure on a space of trajectories, and temporal formulas are subsequently evaluated on those trajectories. The outcome of this evaluation is either a Boolean event or a real-valued functional. While this approach has proven effective in many applications, it maintains a strict separation between the stochastic model and the temporal specification. Temporal properties are not part of the stochastic system itself, but rather act as predicates or functionals defined on its sample paths. As a result, the interaction between temporal logic and probability remains indirect. This separation leads to several limitations. First, temporal properties are not treated as stochastic objects with their own probabilistic structure. Second, logical and temporal operators are not interpreted as operators on stochastic processes, but only as rules for evaluating trajectories. Third, the framework does not naturally support the use of tools from stochastic analysis to study temporal properties, beyond the computation of probabilities of satisfaction. Even in quantitative semantics, where temporal formulas are associated with real-valued robustness functions, these quantities are typically viewed as deterministic mappings applied to trajectories. The resulting values depend on randomness only through the underlying process, rather than being treated as stochastic processes in their own right. 

The above observations and limitations demand a different viewpoint where instead of treating temporal logic as an external specification language, can it be embedded directly into the structure of stochastic processes? In this paper, we propose a new framework, termed \emph{Temporal-Logical Stochastic Processes (TLSP)}, which addresses this question. The central idea is to associate with each temporal logic formula, a stochastic observable derived from the underlying process. In this way, temporal properties are not merely evaluated on trajectories, but are represented as random variables or stochastic processes. This perspective leads to a conceptual shift in how temporal logic interacts with probability. Temporal formulas are interpreted as mappings from trajectories to real numbers, and when applied to a stochastic process, they induce random variables that capture temporal behaviour. Under the proposed novel TLSP framework, it can be observed that logical operations correspond to algebraic operations on these random variables, and temporal operators correspond to transformations acting on stochastic processes. 

The TLSP framework has several important implications as it provides a fully measure-theoretic formulation of temporal logic in stochastic systems, in which temporal properties are embedded within the probabilistic structure. Further, it enables the application of tools from stochastic analysis to the study of temporal properties, including expectation, convergence, and
concentration phenomena. Additionally, the framework reveals an underlying algebraic structure in temporal logic when viewed through the lens of stochastic processes. That is, logical operations such as conjunction and disjunction correspond to operations on random variables, while temporal operators induce structured
transformations of stochastic trajectories. Thus, the proposed TLSP perspective provides a natural bridge between formal methods and probability theory by allowing temporal specifications to be analysed using probabilistic tools, while preserving their logical interpretation. 

The proposed TLSP framework differs fundamentally from existing approaches. In probabilistic model checking, temporal logic defines events whose probabilities are computed \citep{kwiatkowska2007stochastic, baier2008principles}. In stochastic control and reachability analysis, temporal specifications are incorporated into optimization problems
\citep{abate2008probabilistic}. In signal temporal logic, robustness functions are evaluated on deterministic trajectories \citep{donze2010robust, fasoni2020robustness}. In recent years, \citep{Lindemann2021, lindemann2023risk, cleaveland2022risk} and \citep{sleiman_acc} have established probabilistic interpretations of temporal-logic robustness, including robustness distributions, risk measures, and conditions ensuring measurability for discrete- and continuous-time stochastic systems. Building on these developments, the present work systematically studies the entire time-indexed temporal-robustness observable as a stochastic process. Particular emphasis is placed on its pathwise regularity, the topology and convergence of observable processes induced by varying temporal formulas, and the transfer of path-space large-deviation principles to temporal-logical observables. 

The remainder of the paper is organized as follows. Section~\ref{sec_2} introduces preliminaries on stochastic processes and temporal logic. Section~\ref{sec_3} presents the TLSP framework and defines temporal-logical observables. Section~\ref{sec_4} establishes measure theoretical framework for TLSP. Section~\ref{sec_5} discusses the statistical analysis of temporal satisfaction. Section~\ref{sec_6} studies regularity properties of temporal logical observables. Section~\ref{sec_7} presents the distributional properties and extremal characterisation. Section~\ref{sec_8} develops the perturbation studies associated with the robustness of temporal satisfaction. Section~\ref{sec_9} establishes the topology results combined with the convergence of TLOPs.
Section~\ref{sec_10} investigates the large deviation theory results associated with the TLOPs. Finally, the paper is summarized in Section~\ref{sec_11} along with the directions for the future research. 


\section{Notations \& Preliminaries}
\label{sec_2}
\subsection{Preliminaries: Stochastic Processes}
We denote by $\mathbb{R}$ the set of real numbers, by $\mathbb{Z}$ the set of integers, and by $\mathbb{R}_{\geq 0}$ the set of non-negative real numbers. For a measurable space $(E,\mathcal{E})$, the notation $\mathcal{B}(E)$ denotes the Borel $\sigma$-algebra on $E$. For a set $A$, $\mathds{1}_A$ denotes its indicator function. Throughout the manuscript, $(\Omega,\mathcal{F},\mathbb{P})$ denotes a complete probability space. Expectation with respect to $\mathbb{P}$ is denoted by $\mathbb{E}[\cdot]$. For a random variable $Z:\Omega \to \mathbb{R}$, we write $\mathbb{P}(Z \geq 0) := \mathbb{P}(\{\omega \in \Omega : Z(\omega) \geq 0\})$. The time index set is denoted by $\mathbb{T}$, which is either discrete ($\mathbb{Z}_{\geq 0}$) or continuous ($\mathbb{R}_{\geq 0}$). For $t \in \mathbb{T}$ and an interval $I \subseteq \mathbb{T}$, we define the shifted interval as $t + I := \{t + s \mid s \in I\}$. A filtration is a non-decreasing family of sub-$\sigma$-algebras $\{\mathcal{F}_t\}_{t \in \mathbb{T}}$ satisfying $\mathcal{F}_s \subseteq \mathcal{F}_t, \forall s \leq t$. We denote by $\mathcal{C}(\mathbb{T}, \mathcal{X})$ (or $\mathcal{D}(\mathbb{T}, \mathcal{X})$) the space of continuous (or càdlàg) trajectories endowed with the corresponding Borel $\sigma$-algebra. A stochastic process is a collection of measurable maps $X = \{X_t : \Omega \to \mathcal{X}\}_{t \in \mathbb{T}}$, where $(\mathcal{X},\mathcal{B}(\mathcal{X}))$ is a measurable state space. The process is said to be adapted if $X_t$ is $\mathcal{F}_t$-measurable for all $t \in \mathbb{T}$. For $\omega \in \Omega$, the sample path of $X$ is denoted by $X(\omega) := \{X_t(\omega)\}_{t \in \mathbb{T}} \in \mathcal{X}^{\mathbb{T}}$. A functional of trajectories is a measurable mapping $F : \mathcal{X}^{\mathbb{T}} \to \mathbb{R}$. We define the induced random variable $F(X)(\omega) := F(X(\omega))$ and its expectation as 
\begin{align*}
\mathbb{E}[F(X)] := \int_{\Omega} F(X(\omega)) \, d\mathbb{P}(\omega).
\end{align*}
A partial order is a binary relation defining how elements in a set compare, where not all elements need to be comparable. A lattice is a specific type of partially ordered set where any two elements are guaranteed to have both a unique maximum lower bound (meet) and a unique minimum upper bound (join). For two random variables $Z_1, Z_2 \in \mathbb{R}$, the natural pointwise order corresponds to the case when one of the variable is at least as large as the other for every outcome. That is, without loss of generality, 
\begin{align} \label{eqn_partial_order}
Z_1 \leq Z_2 \iff Z_1(\omega) \leq Z_2(\omega), \forall \omega \in \Omega.
\end{align}
Similarly, for two random variables $Z_1, Z_2 \in \mathbb{R}$, they naturally behave like lattices meaning that
\begin{align} \label{eqn_lattices}
\mathrm{meet}(Z_1, Z_2) = \min(Z_1, Z_2), \quad \mathrm{join}(Z_1, Z_2) = \max(Z_1, Z_2). 
\end{align}

\subsection{Temporal Logic Preliminaries}
We use temporal logic as a formal language for describing properties of trajectories of a system in consideration denoted by $\mathbf{\Sigma}$. To familiarise the readers with temporal logic notation, we explain some preliminaries below. A deterministic trajectory is denoted by $x = \{x_t\}_{t\in\mathbb{T}}$ where $x_t\in\mathcal{X}$ is the value of the trajectory at time $t$. On the other hand, for a stochastic process $X=\{X_t\}_{t\in\mathbb{T}}$, each outcome $\omega \in \Omega$ generates one trajectory $X(\omega)=\{X_t(\omega)\}_{t\in\mathbb{T}}$. The symbol $\mu \in \mathcal{A}$ denotes an atomic predicate, meaning a basic statement about the state at one time instant and $\mathcal{A}$ denotes the set of all atomic predicates. Each predicate $\mu \in \mathcal{A}$ is specified by a Borel measurable function $h_\mu:\mathcal{X}\to\mathbb{R}$. The atomic predicate $\mu \in \mathcal{A}$ evaluates to be true at time $t$ along trajectory $x$ if $h_\mu(x_t)\geq 0$. We will denote this satisfaction as $(x,t) \models \mu$ meaning that the trajectory $x$ satisfies the predicate $\mu$ at time $t$. The symbol $\phi$ denotes a temporal formula. The collection of temporal logic formulas is denoted by $\Phi$.
A formula is built recursively from atomic predicates using logical and temporal operators. In this manuscript, we use the grammar
\begin{align} \label{eqn_formula_semantics}
\phi ::= \mu \;|\; \neg\phi \;|\; \phi_1\wedge\phi_2
\;|\; \phi_1 \vee \phi_2 \;|\; \Box_{I}\phi \;|\; \Diamond_{I}\phi .
\end{align}
The notation $(x,t) \models\phi$ means that the trajectory $x$ satisfies the formula $\phi$ at time $t$. The meaning of each notation used in \eqref{eqn_formula_semantics} is as follows: $\mu$ is an atomic predicate; $\neg\phi$ means ``not $\phi$''; $\phi_1 \wedge \phi_2$ means ``$\phi_1$ and $\phi_2$''; $\phi_1 \vee \phi_2 $ means ``$\phi_1$ or $\phi_2$''; $\Box_{I} \phi$ denotes that $\phi$ must \emph{always} hold at all times within the interval $I$ meaning that $\forall s \in t + I, \; (x, s) \models \phi$; $\Diamond_{I} \phi$ denotes that $\phi$ must \emph{eventually} hold at least once within the interval $I$ meaning that $\exists s \in t + I$ such that $(x, s) \models \phi$. Note that we did not use the \emph{until} operator in the above semantics \eqref{eqn_formula_semantics} as it brings in nested optimization based robustness definition which needs a different analysis with a careful inspection and hence is left as a future work. \\

In addition to Boolean satisfaction, we use quantitative semantics called robustness, denoted by
$\rho_\varphi:\mathcal{X}^{\mathbb{T}}\times \mathbb{T} \rightarrow\mathbb{R}$. The quantity $\rho_\varphi(x,t)$ measures the degree to which the trajectory $x$ satisfies or violates the formula $\varphi$ at time $t$. Under the quantitative semantics considered in this manuscript, strictly positive robustness implies satisfaction, while strictly
negative robustness implies violation. That is, 
\begin{align}
\label{eqn_robustness_rho_definition}    
\rho_\phi(x,t) &> 0 \Longrightarrow (x,t) \models \phi, \\
\rho_\phi(x,t) &< 0 \Longrightarrow (x,t) \nvDash \phi.
\end{align}
The converse implications do not necessarily hold. In particular, when $\rho_\varphi(x,t)=0$, the sign of the robustness alone does not determine Boolean satisfaction. For example, consider the atomic predicate $\mu = (x_t\geq 5)$ and a trajectory satisfying $x_t=5$. Then $\rho_\mu(x,t) = 0$, although $(x,t)\models\mu$. On the other hand, for the negated formula $\varphi = \neg\mu$, we have $\rho_\varphi(x,t)=0$,
while $(x,t)\not\models\varphi$. Consequently, zero robustness represents a logical boundary at which Boolean satisfaction must be determined from the
underlying Boolean semantics rather than from the robustness sign alone. This distinction becomes particularly important
when studying satisfaction probabilities, perturbation robustness, and convergence of temporal-logical observables later on in this manuscript. For an atomic predicate $\mu \in \mathcal{A}$ defined by $h_\mu$, we set its robustness as $\rho_\mu(x,t) := h_\mu(x_t)$. The robustness of negation,  conjunction, and disjunction is defined by
\begin{align}
\rho_{\neg\phi}(x,t)
&:= -\rho_\phi(x,t), \label{eqn_robustness_negation}\\
\rho_{\phi_1\wedge\phi_2}(x,t)
&:= \min\{\rho_{\phi_1}(x,t),\rho_{\phi_2}(x,t)\}, \label{eqn_robustness_conjunction} \\
\rho_{\phi_1 \vee \phi_2}(x,t)
&:= \max\{\rho_{\phi_1}(x,t),\rho_{\phi_2}(x,t)\}. \label{eqn_robustness_disjunction}
\end{align}
The robustness of the always operator and the eventually operator over the interval $t+I$ is defined by 
\begin{align}
\rho_{\Box_{I}\phi}(x,t)
&:=
\inf_{s\in t+I}\rho_\phi(x,s), \label{eqn_robustness_always} \\
\rho_{\Diamond_{I}\phi}(x,t)
&:=
\sup_{s\in t+I}\rho_\phi(x,s). \label{eqn_robustness_eventually}
\end{align}
\begin{definition}
The temporal horizon $H_\phi$ of a formula $\phi$ is defined as the maximum amount of future time needed to evaluate a formula $\phi$ at any given time $t \in \mathbb{T}$. 
\end{definition}
The temporal horizon $H_\phi$ for different temporal logical operators are explained in the following Lemma \ref{lemma_temporal_horizons}.
\begin{lemma} \label{lemma_temporal_horizons}
For a temporal logic formula $\phi$ generated by \eqref{eqn_formula_semantics}, assume that every temporal interval $I$ is bounded and define $\overline I := \sup I$. Then, the temporal horizon $H_\phi$ is given recursively as 
\begin{align}
H_{\mu} &= 0, \label{eqn_temporal_horizon_predicate} \\
H_{\neg \phi} &= H_\phi, \label{eqn_temporal_horizon_negation} \\
H_{\phi_1\wedge\phi_2} &= H_{\phi_1\vee\phi_2} = \max\{H_{\phi_1},H_{\phi_2}\}, \label{eqn_temporal_horizon_conjunction_disjunction} \\
H_{\Box_{I}\phi} &= H_{\Diamond_{I}\phi} = \overline I + H_\phi. \label{eqn_temporal_horizon_always_eventually}
\end{align}
\end{lemma}
\begin{proof}
For an atomic predicate $\mu$, satisfaction at time $t$ depends only on the state at time $t$. Hence no future information is needed, and $H_\mu = 0$. For negation, the formula $\neg\phi$ is evaluated by first evaluating $\phi$ at the same time and then reversing its truth value, or equivalently changing the sign of its robustness value. Therefore, negation does not introduce any additional time dependence, and $H_{\neg\phi} = H_\phi$. Both for conjunction and disjunction, the formula $\phi_1 \wedge \phi_2$ and $\phi_1\vee\phi_2$ respectively are evaluated by evaluating both $\phi_1$ and $\phi_2$ at the same time. Thus, the amount of future information required is the larger of the two horizons. Hence,
$H_{\phi_1\wedge\phi_2} = H_{\phi_1\vee\phi_2} =
\max\{H_{\phi_1},H_{\phi_2}\}$. Now consider the always formula $\Box_{I}\phi$. To evaluate
$\Box_{I}\phi$ at time $t$, one must evaluate $\phi$ at every time $s \in t + I$. Since $s$ can be as large as $t + \overline I$, and evaluating
$\phi$ at such an $s$ may require information up to time $s + H_\phi$, the largest time needed is $t + \overline I + H_\phi$. Therefore, relative to the initial time $t$, the horizon is $H_{\Box_{I}\phi} =
\overline I + H_\phi$. Similarly, the eventually formula $\Diamond_{I}\phi$ is evaluated by checking whether $\phi$ holds at least once for some $s \in t + I$. Since the possible time $s$ can again be as large as $t + \overline I$, evaluating the formula may require information up to time $t+\overline I+H_\phi$. Hence, relative to the initial time $t$, the horizon is $H_{\Diamond_{I}\phi} =
\overline I+H_\phi$. 
\end{proof}
\noindent Using Lemma \ref{lemma_temporal_horizons}, we can evaluate the temporal horizon of any temporal formula $\phi \in \Phi$. For example, $H_{\mu} = 0, H_{\Box_{[0,T]}\mu} = T$, and $H_{\Diamond_{[0,T]}(\Box_{[0,S]}\mu)} =  T + S$. 

\section{Temporal-Logical Stochastic Processes (TLSP)}\label{sec_3}
Classically, temporal logic is evaluated after a stochastic process generates a trajectory. That is, for a stochastic process $X$, one typically studies quantities of the form $\mathbb{P}(X \models \phi)$ and in this viewpoint the temporal logic is treated as an external specification. In contrast, we seek to develop a framework in which temporal logic is embedded 
directly into the stochastic structure. To this end, we first define the temporal logical observables and the associated Temporal-Logical Observable Process (TLOP) below.

\subsection{TLOP and TLSP}
We now propose a change from the external specification viewpoint (evaluating a stochastic process against a formula specified using temporal logics) by taking the deterministic robustness functional $\rho_{\phi}(x,t)$ and replacing the deterministic trajectory $x$ by the random trajectory $X(\omega)$. Since, the trajectory $X(\omega)$ is random, the resulting robustness functional $\rho_{\phi}(X(\omega),t)$ will also be random as well as it depends upon the outcome $\omega$. The robustness function $\rho_{\phi}$ through its magnitude capturing the margin of satisfaction or violation of a formula $\phi$ has turned random with this change. The formal definition is given below.     
\begin{definition}
For each temporal formula $\phi\in\Phi$ and each time $t\in\mathbb{T}$, define the temporal-logical observable mapping induced by the stochastic process $X = \{X_t\}_{t\in\mathbb{T}}$ and the formula $\phi$ as $S_{\phi}^{X}(t): \Omega \to \mathbb{R}$ and it is given by
\begin{align}
    S_{\phi}^{X}(t)(\omega)
    :=
    \rho_{\phi}\big(X(\omega),t\big),
    \qquad \omega\in\Omega.
    \label{eqn_temporal_logical_observable}
\end{align}
\end{definition}
\noindent Note that, $X(\omega) = \{X_s(\omega)\}_{s\in\mathbb{T}}$ denotes the sample path of the stochastic process $X$. Under the measurability conditions established later in Theorem \ref{thm_TLOP_measurability} of this manuscript, $S_{\phi}^{X}(t)$ is a real-valued random variable for every $t \in \mathbb{T}$. Its value $S_{\phi}^{X}(t)(\omega)$ quantifies the temporal robustness of the sample path $X(\omega)$ with respect to the formula $\phi$ at time $t$.
\begin{definition}
Suppose that $S_{\phi}^{X}(t)$ is measurable $\forall t\in\mathbb{T}$. Then, the stochastic process $S_\phi^X := \left\{S_\phi^X(t)\right\}_{t\in\mathbb{T}}$ is called the TLOP induced by the process $X$ and formula $\phi$.
\end{definition}
\noindent Related stochastic interpretations of temporal-logic robustness have appeared previously. Bartocci et al.~\citep{Bartocci2015} extended quantitative robustness to stochastic models and studied the resulting distribution of robustness degrees and their average values. More recently, Lindemann et al. in ~\citep{Lindemann2021} considered STL robustness over discrete-time stochastic processes and quantified the associated robustness risk using general risk measures. These works establish that temporal-logic robustness can itself be treated as a random quantity. The perspective pursued here is complementary: rather than focusing only on the robustness random variable associated with a fixed evaluation time, we regard the entire time-indexed family
$\{S_{\phi}^{X}(t)\}_{t\in\mathbb{T}}$ as a stochastic process and study
its measurability, pathwise regularity, perturbation properties, and process-level convergence. It is clear from the above definition that while $X = \{X_t\}_{t \in \mathbb{T}}$ denotes the original physical process, the notation $S_\phi^X := \left\{S_\phi^X(t)\right\}_{t\in\mathbb{T}}$ denotes the new stochastic process describing the temporal-logical behaviour of $X$ with respect to the formula $\phi$. To understand this better, let $X_t$ denote the temperature process of some chemical plant and assume that the process $X_t$ needs to satisfy a simple temporal specification that ``from time $t$, the temperature must remain at least $20^\circ$ for the next six units of time''. Mathematically, we can write this specification $\phi$ using the temporal logical operators in \eqref{eqn_formula_semantics} as
\begin{align*}
\phi
=
\Box_{[0,5]} (X \geq 20^\circ)
\implies
S_\phi^X(t)
=
\inf_{s \in [t, t+5]} (X_s - 20). 
\end{align*}
Consider the following two future realisations of the temperature process at time $t$ namely
\begin{align*}
X(\omega_1) 
&= \{21,23,22,21,24\} \implies S_\phi^X(t) = 1, \text{ and }\\
X(\omega_2) 
&= \{21,23,19,22,25\} \implies S_\phi^X(t) = -1
\end{align*}
respectively. That is, different outcomes $\omega_1, \omega_2$ result in different values for $S_\phi^X(t)$. So, $S_\phi^X(t)$ will have a probability distribution. Then, one can study quantities like moments of $S_\phi^X(t)$, $\mathbb{P}(S_\phi^X(t) > c)$ for some $c \in \mathbb{R}$ and so on. 
\begin{lemma}
\label{lemma_probability_formula_satisfaction}
Let $X=\{X_t\}_{t\in\mathbb{T}}$ be a stochastic process and let
$\phi\in\Phi$. For every $t\in\mathbb{T}$,
\begin{align}
    \{S_\phi^X(t)>0\}
    \subseteq
    \{\omega\in\Omega \mid X(\omega),t \models\phi\}
    \subseteq
    \{S_\phi^X(t)\geq 0\} 
    \implies     \mathbb{P}\big(S_\phi^X(t)>0\big)
    \leq
    \mathbb{P}\big(X,t \models\phi\big)
    \leq    \mathbb{P}\big(S_\phi^X(t)\geq0\big).
    \label{eqn_probability_formula_satisfaction_bounds}
\end{align}
In particular, if $\mathbb{P}\big(S_\phi^X(t)=0\big)=0$, then 
\begin{align}
\label{eqn_probability_formula_satisfaction}
    \mathbb{P}\big(X,t \models\phi\big)
    =
    \mathbb{P}\big(S_\phi^X(t)>0\big)
    =
    \mathbb{P}\big(S_\phi^X(t)\geq0\big).
\end{align}
\end{lemma}
\begin{proof}
We know that positive robustness implies
satisfaction, while negative robustness implies violation. So,
\begin{align*}
\{S_\phi^X(t)>0\}
\subseteq
\{\omega:(X(\omega),t)\models\phi\} \\
\{S_\phi^X(t)<0\}
\subseteq
\{\omega:(X(\omega),t) \nvDash \phi\}
\end{align*}
Hence, $\{S_\phi^X(t)>0\} \subseteq \{\omega \mid (X(\omega),t) \models\phi\} \subseteq \{S_\phi^X(t)\geq0\}$. Taking probabilities gives \eqref{eqn_probability_formula_satisfaction_bounds}. If
$\mathbb{P}(S_\phi^X(t)=0)=0$, then the events
$\{S_\phi^X(t)>0\}$ and $\{S_\phi^X(t)\geq0\}$ differ only by a
null event, and the equality in \eqref{eqn_probability_formula_satisfaction} follows.
\end{proof}
Note that the probability of satisfaction can also be inferred using the indicator function approach. To this end, we define the indicator function associated with the non-negativity of TLOP below.
\begin{definition}
For each $t\in\mathbb{T}$ and temporal formula $\phi\in\Phi$, define the
temporal satisfaction indicator associated with the stochastic process $X$ as
\begin{align}
\label{eqn_indicator_formula_satisfaction_definition}
    \mathds{1}_{\phi}^{X}(t)
    :=
    \mathds{1}_{
        \{\omega\in\Omega:
        (X(\omega),t)\models\phi\}
    } 
    \iff 
    \mathds{1}_{\phi}^{X}(t)(\omega)
    =
    \begin{cases}
        1,
        & \text{if } (X(\omega),t)\models\phi,\\
        0,
        & \text{otherwise}.
    \end{cases},
    \quad 
    \forall \omega\in\Omega.
\end{align}
\end{definition}
\noindent Thus, $\mathds{1}_{\phi}^X(t)$ is a Bernoulli random variable that records if the stochastic trajectory satisfies the formula $\phi$ at time $t$. 
\begin{corollary} \label{corollary_indicator_formula_satisfaction}
Using \eqref{eqn_probability_formula_satisfaction} and \eqref{eqn_indicator_formula_satisfaction_definition}, we see that, $\mathbb{P}(X, t \models \phi) = \mathbb{E}\big[\mathds{1}_{\phi}^X(t)\big]$. 
\end{corollary}
\begin{proof}
The proof is straightforward and hence is omitted.
\end{proof}
Through the above development, we saw that the satisfaction probability of a temporal logic specification by a stochastic process has become trivial when expressed in terms of the TLOP language.
\begin{definition}
Given a collection of temporal logic formulas $\Phi$ and an $\mathcal{X}$-valued stochastic process $X=\{X_t\}_{t\in\mathbb{T}}$ defined on $(\Omega,\mathcal F,\mathbb{P})$, a Temporal-Logical Stochastic Process (TLSP) is denoted by the pair $(X,\Phi)$.
\end{definition}
\noindent For each formula $\phi\in\Phi$, the TLSP induces the TLOP denoted by $S_\phi^X = \left\{S_\phi^X(t)\right\}_{t\in\mathbb{T}}$. Thus, a TLSP is not merely a stochastic process together with external
specifications, rather it is a stochastic process equipped with a family of
measurable temporal-logical observables.


\section*{Illustrative Example: A TLSP for Stock Prices}
To illustrate the proposed TLSP framework, consider the evolution of the price of a
single financial asset, such as a stock, over a finite investment horizon. Let
$S_t$ denote the price of the stock at time $t \in \mathbb{T}_{\geq 0}$. Since the stock price is
uncertain and evolves randomly over time, it is natural to model
$\{S_t\}_{t\geq 0}$ as a stochastic process. A classical model for such an
evolution is the geometric Brownian motion model, in which the stock price
satisfies the stochastic differential equation
\begin{align}
dS_t = \mu S_t\,dt + \sigma S_t\,dW_t.
\end{align}
Here, $\mu\in\mathbb R$ denotes the drift coefficient representing the average instantaneous rate of growth of the stock price. The parameter $\sigma > 0$ denotes the volatility coefficient, which measures the intensity of random fluctuations in the stock price. The process $\{W_t\}_{t\geq 0}$ is a standard Brownian motion \citep{uhlenbeck1930theory} defined on a probability space $(\Omega, \mathcal F,\mathbb{P})$. Each outcome $\omega \in \Omega$ corresponds to one possible realization of market randomness, and therefore generates one possible stock-price trajectory $S(\omega) = \{ S_t(\omega)\}_{t\geq 0}$. Thus, for each fixed time $t \in \mathbb{T}_{\geq 0}$, the proice of stock $S_t$ is a random variable, while for each fixed outcome $\omega$, $S_t(\omega)$ describes one possible price path. In the notation of the TLSP framework, we take the underlying stochastic process to be $X_t := S_t$, so that $X=\{X_t\}_{t\geq 0}$ represents the stock-price stochastic process. \\

We now define a temporal specification that expresses a meaningful investment requirement. Suppose that an investor is interested in the following stock evolution property that over the next $T > 0$ units of time, the stock price should eventually rise above a desired target level $K > 0$, while never falling below an unacceptable safety threshold level $L > 0$. Implicitly, one should interpret this in the natural case in which the desired target price is above the minimum acceptable price meaning that $K > L$. This investment requirement contains two distinct temporal components. The first component is an eventual-gain requirement where the the stock should reach or
exceed the target level $K$ at least once during the interval $[t,t+T]$. The second component is a safety requirement in which the stock should remain above the lower threshold $L$ throughout the same interval. Using the notations explained previously in this paper, the two temporal specifications will give rise to the following temporal formula given by
\begin{align}
\label{eqn_stock_example_formula}
\phi
:=
\Diamond_{[0,T]}\mu_K
\;\wedge\;
\Box_{[0,T]}\mu_L,
\end{align}
where the atomic predicate $\mu_K$ represents the statement that the stock price is at least $K$, and the atomic predicate $\mu_L$ represents the statement that the stock price is at least $L$. These predicates are defined through the functions
\begin{align}
h_{\mu_K}(x) := x-K,
\qquad
h_{\mu_L}(x) := x-L.
\end{align}
Therefore, at a time $s$, the predicate $\mu_K$ is satisfied if $S_s - K \geq 0$, which is equivalent to $S_s \geq K$. Similarly, the predicate $\mu_L$ is satisfied if $S_s - L \geq 0$, which is equivalent to $S_s \geq L$. The formula $\phi$ in \eqref{eqn_stock_example_formula} states that the stock must reach the target level $K$ at least once during the future interval $[0,T]$, and must remain above the safety threshold $L$ throughout that interval $[0,T]$. We now construct the temporal-logical observable associated with this formula.
For a deterministic price trajectory $x=\{x_s\}_{s\geq 0}$, the robustness of
the eventual-gain requirement at time $t$ and the robustness of the safety requirement are respectively given by
\begin{align}
\rho_{\Diamond_{[0,T]}\mu_K}(x,t)
&=
\sup_{s\in[t,t+T]}(x_s - K), \label{eqn_stock_gain_spec_robustness}\\
\rho_{\Box_{[0,T]}\mu_L}(x,t)
&=
\inf_{s\in[t,t+T]}(x_s - L). \label{eqn_stock_safety_spec_robustness}
\end{align}
The quantity $\rho_{\Diamond_{[0,T]}\mu_K}(x,t)$ is positive precisely when the price trajectory exceeds the target level $K$ at least once during the interval $[t,t+T]$. The larger this quantity is, the more strongly the target-reaching requirement is satisfied. Similarly, the quantity $\rho_{\Box_{[0,T]}\mu_L}(x,t)$ is positive precisely when the trajectory remains above the
threshold $L$ for all times in the interval $[t,t+T]$. The infimum appears because the always requirement is governed by the worst value attained over the interval. Since both the temporal specifications are required to be satisfied, we combine the two individual temporal robustness \eqref{eqn_stock_gain_spec_robustness} and \eqref{eqn_stock_safety_spec_robustness} using conjunction to get the robustness of the formula $\phi$ in \eqref{eqn_stock_example_formula} as
\begin{align}
\label{eqn_stock_formula_robustness}
\rho_\phi(x,t)
=
\min\left\{
\sup_{s\in[t,t+T]}(x_s-K),
\inf_{s\in[t,t+T]}(x_s-L)
\right\}.
\end{align}
The minimum in \eqref{eqn_stock_formula_robustness} appears because the conjunction is satisfied only as strongly as its weaker component. Applying this robustness functional to the stochastic stock-price process produces the temporal logical observable
\begin{align}
S_\phi^X(t)(\omega)
=
\min\left\{
\sup_{s\in[t,t+T]}(S_s(\omega)-K),
\inf_{s\in[t,t+T]}(S_s(\omega)-L)
\right\}.
\end{align}
For each fixed time $t$, the quantity $S_{\phi}^{X}(t)$ is a random variable. It assigns to each possible market outcome $\omega$ a real number describing how strongly the corresponding stock-price trajectory
satisfies or violates the temporal investment requirement. For the particular formula $\phi$ in \eqref{eqn_stock_example_formula}, which is constructed only
from threshold predicates, conjunction, and the temporal operators
$\lozenge$ and $\square$, Boolean satisfaction can be recovered exactly
from the sign of the TLOP. Indeed,
\begin{align*}
S_{\phi}^{X}(t)(\omega)\geq 0
&\Longleftrightarrow
\sup_{s\in[t,t+T]}
\left(S_s(\omega)-K\right)\geq 0
\quad\text{and}\quad
\inf_{s\in[t,t+T]}
\left(S_s(\omega)-L\right)\geq 0
\\
&\Longleftrightarrow
\exists s\in[t,t+T]\text{ such that }S_s(\omega)\geq K,
\quad\text{and}\quad
S_s(\omega)\geq L,\ \forall s\in[t,t+T]
\\
&\Longleftrightarrow
X(\omega),t \models\phi.
\end{align*}
Hence, $\mathbb{P}\big(X,t \models\phi\big) = \mathbb{P}\big(S_{\phi}^{X}(t)\geq0\big)$. This probability has a natural financial interpretation in the sense that it is the probability that the stock satisfies the desired temporal investment profile over the future time window $[t,t+T]$. Importantly, TLSP does not merely ask whether the formula satisfaction event occurs. It also constructs a stochastic observable whose distribution contains more information than the satisfaction probability alone. For example, the expectation of $S_\phi^X(t)$ describes the average robustness of the temporal specification, while the lower tail of its distribution describes the risk of severe violation. The indicator random variable associated with temporal satisfaction is $\mathds{1}_\phi^X(t) =
\mathds{1}_{\{S_\phi^X(t) \geq 0\}}$. This variable takes the value one when the temporal specification is satisfied and zero otherwise.  This example demonstrates the central idea of TLSP in which the temporal formula $\phi$ given by \eqref{eqn_stock_example_formula} is not treated merely as an external statement about the stock-price trajectory. Instead, it induces a new stochastic process $\{S_\phi^X(t)\}_{t\geq 0}$ whose values quantify temporal satisfaction over
moving future windows. In this sense, the temporal specification becomes an intrinsic stochastic observable associated with the original price process. The same construction can be applied to other stochastic systems such as biological processes, climate models, energy systems, and multi-agent systems, whenever one wishes to study their temporally structured properties under uncertainty.

\section{Measure Theoretic Framework for TLSP} \label{sec_4}
While the manuscript aims at elevating the temporal logic from an external verification tool to an intrinsic stochastic observable, the central objective of this section is to enable a unified measure-theoretic treatment of temporal properties of stochastic systems. While a preliminary work in this direction already exists in \citep{lindemann2023risk}, we develop a measure-theoretic framework for TLSP in which temporal logic formulas induce well-defined stochastic observables with the following properties:
\begin{enumerate}
    \item For each temporal formula $\phi$ and time $t\in\mathbb{T}$, the mapping $S_\phi^X(t):\Omega \to \mathbb{R}$ is a measurable random variable. If $X$ is adapted to a filtration $\{\mathcal F_t\}$, then the induced observable process $\{S_\phi^X(t)\}$ is adapted to a shifted filtration determined by $H_\phi$. Subsequently, the temporal properties of the stochastic system can be analysed via probabilistic quantities such as $\mathbb{E}[S_\phi^X(t)]$, $\mathbb{P}(S_\phi^X(t)\geq 0)$ enabling the use of tools from stochastic analysis.
    \item The temporal operators $\Box_{I}$ and $\Diamond_{I}$ correspond to measurable transformations of stochastic processes.
    \item The class of temporal-logical observables is closed under Boolean operations and analogous relations hold for negation.
\end{enumerate}
\subsection{Measurability of TLOP}
Consider the formula $\phi=\Diamond_{[0,10]}(X \geq 5)$. To determine the $S^X_{\phi}(t)$, we need $X_s$ for $t \leq s \leq t + 10$. Though $S^X_{\phi}(t)$ is expressed at time $t$, we need the future of $X$ to determine $S^X_{\phi}(t)$ meaning that $S^X_{\phi}(t)$ is generally not measurable with respect to the filtration $\mathcal{F}_t$. Using the temporal horizon, we will now prove the measurability of the temporal logical stochastic observables with respect to their corresponding shifted filtration in the following theorem, tailored for discrete-time setting. 
\begin{theorem}
\label{thm_TLOP_measurability}
Let $\mathbb{T}=\mathbb{Z}_{\geq 0}$ and let
$X=\{X_t\}_{t\in\mathbb{T}}$ be adapted to the filtration $\{\mathcal{F}_t\}_{t\in\mathbb{T}}$. Assume that every temporal interval $I$ appearing in the formula $\phi\in\Phi$ is finite. Then, $S_\phi^X(t)$ is $\mathcal{F}_{t+H_\phi}$-measurable $\forall t\in\mathbb{T}$.
\end{theorem}
\begin{proof}
We will prove the result by structural induction on elements constituting the formula $\phi$ given by \eqref{eqn_formula_semantics}. 
\begin{enumerate}
    \item Case 1: Consider the case when the formula is specified as an atomic predicate meaning that $\phi = \mu$. Then $S_\mu^X(t)(\omega) = \rho_\mu(X(\omega),t) = h_\mu(X_t(\omega))$. Since $X$ is adapted, $X_t$ is $\mathcal F_t$-measurable. Since $h_\mu$ is Borel measurable, $h_\mu(X_t)$ is $\mathcal F_t$-measurable. Since $H_\mu = 0$ from \eqref{eqn_temporal_horizon_predicate} in Lemma \ref{lemma_temporal_horizons}, the claim follows.
    \item Case 2: Consider the case when the formula is specified as a negation meaning that $\phi = \neg \psi$. Then $S_{\neg \psi}^X(t) = - S_\psi^X(t)$. By the induction hypothesis, $S_\psi^X(t)$ is $\mathcal F_{t+H_\psi}$-measurable. Hence $S_{\neg \psi}^X(t)$ is also $\mathcal F_{t+H_\psi}$-measurable. Since $H_{\neg \psi}=H_\psi$ from \eqref{eqn_temporal_horizon_negation} in Lemma \ref{lemma_temporal_horizons}, the claim follows.
    \item Case 3: Consider the case when the formula is specified as a conjunction or as a disjunction. Note that $H_{\psi_1 \wedge \psi_2} = H_{\psi_1 \vee \psi_2} = \max\{H_{\psi_1},H_{\psi_2}\}$ from \eqref{eqn_temporal_horizon_conjunction_disjunction} in Lemma \ref{lemma_temporal_horizons}. Let $\phi_1 = \psi_1 \wedge \psi_2$ and $\phi_2 = \psi_1 \vee \psi_2$. Then
    \begin{align}
    S_{\phi_1}^X(t)
    &=
    S_{\psi_1 \wedge \psi_2}^X(t)
    =
    \min\{S_{\psi_1}^X(t), S_{\psi_2}^X(t)\} \\
    S_{\phi_2}^X(t)
    &=
    S_{\psi_1 \vee \psi_2}^X(t)
    =
    \max\{S_{\psi_1}^X(t), S_{\psi_2}^X(t)\}.
    \end{align}
    By the induction hypothesis, $S_{\psi_1}^{X}(t)$ is $\mathcal{F}_{t+H_{\psi_1}}$-measurable and $S_{\psi_2}^{X}(t)$ is $\mathcal{F}_{t+H_{\psi_2}}$-measurable. Since the filtration is non-decreasing,
    \begin{align*}
    \mathcal{F}_{t+H_{\psi_j}}
    \subseteq
    \mathcal{F}_{t+\max\{H_{\psi_1},H_{\psi_2}\}}, 
    \quad \text{for } j = \{1, 2\}.
    \end{align*}
    Hence both $S_{\psi_1}^{X}(t)$ and $S_{\psi_2}^{X}(t)$ are $\mathcal{F}_{t+\max\{H_{\psi_1},H_{\psi_2}\}}$-measurable. Since the minimum and maximum mappings are Borel measurable, the result follows.
    \item Case 4: Consider the case when the formula is specified using always operator meaning that $\phi = \Box_{I} \psi$. Then $S_{\Box_{I} \psi}^X(t) = \min_{s\in t+I} S_\psi^X(s)$. For each $s\in t+I$, the induction hypothesis implies that $S_\psi^X(s)$ is $\mathcal F_{s+H_\psi}$-measurable. Since $s \leq t + \max I$, we have, $s+H_\psi \leq t + \max I + H_\psi$ from Lemma \ref{lemma_temporal_horizons}. Thus, each $S_\psi^X(s)$ is measurable with respect to $\mathcal F_{t+\max I+H_\psi}$. The infimum of a finite collection of measurable random variables is measurable. Hence $S_{\Box_{I} \psi}^X(t)$ is $\mathcal F_{t+H_{\Box_{I} \psi}}$-measurable.
    \item Case 5: Consider the case when the formula is specified using eventually operator meaning that $\phi = \Diamond_{I} \psi$. The proof is analogous to the always operator case. Then $S_{\Diamond_{I} \psi}^X(t) = \max_{s\in t+I} S_\psi^X(s)$. Using the same argument of always operator above from Lemma \ref{lemma_temporal_horizons}, we conclude that $S_{\Diamond_{I} \psi}^X(t)$ is $\mathcal F_{t+H_{\Diamond_{I} \psi}}$-measurable.
\end{enumerate}
Since every formula $\phi \in \Phi$ is obtained by finitely many applications of the above five cases, the proof is complete.
\end{proof}
\begin{remark}
The continuous-time case requires additional regularity assumptions. In particular, when $\mathbb{T}=\mathbb{R}_{\geq 0}$ and the temporal intervals are compact, measurability of the temporal infimum and supremum can be established on suitable path spaces, such as
$D(\mathbb{R}_{\geq 0},\mathbb{R})$, equipped with its Borel $\sigma$-algebra. This setting is treated separately in the next subsection. 
\end{remark}
\subsection{Temporal operators as measurable operators on process space}
\begin{theorem} \label{thm_temporal_operators_measurable}
Let $\mathbb{T}$ be either $\mathbb{Z}_{\geq 0}$ or $\mathbb{R}_{\geq 0}$, and let $I$ be a
non-empty temporal interval. In discrete time, let $\mathcal{Z} = \mathbb{R}^{\mathbb{Z}_{\geq 0}}$
be equipped with the product $\sigma$-algebra, and assume that $I$ is finite.
In continuous time, let $\mathcal{Z} = D(\mathbb{R}_{\geq 0},\mathbb{R})$ be the space of real-valued càdlàg trajectories, equipped with the Skorokhod $J_1$ topology and its Borel $\sigma$-algebra. Assume that
$I=[a,b]$, where $0\leq a\leq b<\infty$.
Define the temporal operators by
\begin{align}
(\mathcal{T}_I^\Box z)(t)
&:=
\inf_{s\in t+I}z(s), \\
(\mathcal{T}_I^\Diamond z)(t)
&:=
\sup_{s\in t+I}z(s).
\end{align}
Then the following statements hold.
\begin{enumerate}
\item In discrete time, both operators are measurable mappings from $\mathcal{Z}$ into $\mathbb{R}^{\mathbb{T}}$ equipped with the product $\sigma$-algebra.
\item In continuous time, both operators map càdlàg trajectories into càdlàg trajectories. Furthermore, $\mathcal{T}_I^\Box,\mathcal{T}_I^\Diamond: \mathcal{Z}\longrightarrow\mathcal{Z}$ are Borel measurable mappings with respect to the Skorokhod $J_1$ Borel $\sigma$-algebra.
\end{enumerate}
Moreover, for every temporal formula $\varphi$ whose induced TLOP $S_\varphi^X$ has sample paths in $\mathcal{Z}$, the following identities hold pathwise wherever the corresponding
sample-path assumptions are satisfied:
\begin{equation}
S_{\Box_I\varphi}^X
=
\mathcal{T}_I^\Box(S_\varphi^X),
\qquad
S_{\Diamond_I\varphi}^X
=
\mathcal{T}_I^\Diamond(S_\varphi^X).
\end{equation}
\end{theorem}

\begin{proof}
We first establish the discrete-time result. For any fixed $t\in\mathbb{Z}_{\geq 0}$, the set $t+I$ is finite. Therefore,
\begin{align*}
(\mathcal{T}_I^\Box z)(t)
&=\min_{s\in t+I}z(s),\\
(\mathcal{T}_I^\Diamond z)(t)
&=\max_{s\in t+I}z(s).
\end{align*}
Since coordinate evaluations are measurable under the product $\sigma$-algebra, and finite minima and maxima of measurable functions are
measurable, both temporal operators are
coordinate-wise measurable. Hence they are measurable into
$\mathbb{R}^{\mathbb{T}}$ equipped with the product $\sigma$-algebra.

Now consider the continuous time setting. Fix $z\in D(\mathbb{R}_{\geq 0},\mathbb{R})$
and a finite observation horizon $R>0$.
A càdlàg function on a compact interval can be uniformly approximated by right-continuous step functions having finitely many jumps \citep[Theorem~2.3]{Frankova2019}. Consequently, there exists a sequence
$\{z_n\}_{n\in\mathbb{N}}$ of such step
functions satisfying
\begin{equation*}
\sup_{s\in[0,R+b]}
|z_n(s)-z(s)|
\longrightarrow 0.
\end{equation*}
For every pair of bounded functions $f,g$, the elementary inequalities
\begin{align*}
\left|
\inf_{s\in t+I}f(s)
-
\inf_{s\in t+I}g(s)
\right|
&\leq
\sup_{s\in t+I}|f(s)-g(s)|,\\
\left|
\sup_{s\in t+I}f(s)
-
\sup_{s\in t+I}g(s)
\right|
&\leq
\sup_{s\in t+I}|f(s)-g(s)|
\end{align*}
hold. So, for either temporal operator
$\mathcal{T}_I$,
\begin{equation*}
\sup_{t\in[0,R]}
\left|
(\mathcal{T}_I z_n)(t)
-
(\mathcal{T}_I z)(t)
\right|
\leq
\sup_{s\in[0,R+b]}|z_n(s)-z(s)|.
\end{equation*}
Thus, the transformed step functions converge uniformly to $\mathcal{T}_I z$ on $[0,R]$. Since each $z_n$ is a right-continuous step function with finitely many jumps, its moving infimum and supremum over the fixed interval
$I=[a,b]$ are again right-continuous step functions with finitely many jumps.
Hence $\mathcal{T}_I z_n$ is càdlàg. Since a uniform limit of càdlàg functions on a compact interval is càdlàg, $\mathcal{T}_I z$ is càdlàg on $[0,R]$. As $R>0$ was arbitrary, both
$\mathcal{T}_I^\Box z$ and
$\mathcal{T}_I^\Diamond z$ belong to
$D(\mathbb{R}_{\geq0},\mathbb{R})$. \\

\noindent Fix $t\geq0$. Since $z$ is right-continuous, we have
\begin{align*}
\inf_{s\in[t+a,t+b]}z(s)
&=
\inf_{s\in Q_{t,I}}z(s),\\
\sup_{s\in[t+a,t+b]}z(s)
&=
\sup_{s\in Q_{t,I}}z(s),
\end{align*}
where
\begin{equation*}
Q_{t,I}
:=
\left(
\mathbb{Q}\cap[t+a,t+b)
\right)
\cup\{t+b\}.
\end{equation*}
The set $Q_{t,I}$ is countable. For every fixed $s\geq0$, the coordinate
evaluation $z\longmapsto z(s)$ is Borel measurable on the Skorokhod space
$D(\mathbb{R}_{\geq0},\mathbb{R})$. Since countable infima and suprema of
measurable real-valued functions are
measurable, the mappings $z\longmapsto(\mathcal{T}_I^\Box z)(t)$ 
and $z\longmapsto(\mathcal{T}_I^\Diamond z)(t)$ are Borel measurable for every fixed $t$. The Borel $\sigma$-algebra of the Skorokhod space is generated by its coordinate evaluation mappings at a countable dense set of time points. Since we already established that both temporal operators take values in the Skorokhod space, their coordinate-wise measurability implies that $\mathcal{T}_I^\Box,\mathcal{T}_I^\Diamond: \mathcal{Z}\longrightarrow\mathcal{Z}$ are Borel measurable. \\

\noindent By the quantitative semantics of the always and eventually operators, we see that 
\begin{align*}
S_{\Box_I\varphi}^X(t)
&=
\inf_{s\in t+I}S_\varphi^X(s)
=
\left(
\mathcal{T}_I^\Box S_\varphi^X
\right)(t), \\
S_{\Diamond_I\varphi}^X(t)
&=
\sup_{s\in t+I}S_\varphi^X(s)
=
\left(
\mathcal{T}_I^\Diamond S_\varphi^X
\right)(t).
\end{align*}
Since these identities hold $\forall t$, the result follows and the proof is complete.
\end{proof}

The continuous-time result requires that the process to which a temporal operator is applied has càdlàg sample paths. In particular, Borel measurability of the atomic predicate functions alone does not guarantee càdlàg temporal-logical observables. A sufficient condition is that the underlying process $X$ has càdlàg sample paths almost surely and that all atomic predicate functions $h_\mu$ are continuous. Under these assumptions, the atomic robustness processes $h_\mu(X_t)$ are càdlàg. Negation, conjunction, and disjunction preserve
càdlàg paths, while Theorem \ref{thm_temporal_operators_measurable} establishes the same property for the temporal operators. Consequently, structural induction shows that every TLOP generated by the formula grammar
has càdlàg sample paths almost surely.
Moreover, if $X$ is a measurable random element of the Skorokhod path space, the corresponding TLOP is also a measurable random element of
that space. If $X$ is adapted to a filtration $\{\mathcal{F}_t\}_{t\geq0}$, then for every formula $\varphi$ with finite temporal horizon, $S_\varphi^X(t)$ is $\mathcal{F}_{t+H_\varphi}$-measurable. 
The latter property follows by structural induction, using the countable representation of the temporal infimum and supremum
established in the proof of Theorem \ref{thm_temporal_operators_measurable}.


\subsection{Algebra of Temporal-Logical Observables}
We will now study the structural properties of the mapping $\phi \mapsto S_\phi^X(t)$ which associates temporal logic formulas with real-valued random variables.

\subsubsection{Pointwise Algebraic Structure}

Let $\phi_1$ and $\phi_2$ be temporal logic formulas. The following result
characterizes the algebraic closure explaining how logical operations translate into operations on temporal-logical observables.
\begin{proposition} \label{proposition_ptwise_algebraic_structure}
For every $t\in\mathbb{T}$, the temporal-logical observables satisfy
\begin{align}
S_{\phi_1 \wedge \phi_2}^X(t)
&=
\min\{S_{\phi_1}^X(t), S_{\phi_2}^X(t)\},\\
S_{\phi_1 \vee \phi_2}^X(t)
&=
\max\{S_{\phi_1}^X(t), S_{\phi_2}^X(t)\},\\
S_{\neg \phi}^X(t)
&=
- S_\phi^X(t).
\end{align}
\end{proposition}

\begin{proof}
The result follows directly by substituting $x=X(\omega)$ in the definition of the robustness semantics in \eqref{eqn_robustness_negation}-\eqref{eqn_robustness_disjunction}.
\end{proof}

\subsubsection{Lattice Structure}
\noindent The previous proposition implies that the collection of temporal-logical
observables inherits a natural lattice structure.
\begin{proposition}
Fix $t\in\mathbb{T}$. The collection
$\mathcal S_t^X := \left\{S_\phi^X(t) \mid \phi \text{ temporal formula}\right\}$ is closed under the pointwise operations given by \eqref{eqn_lattices}. In particular, $\mathcal S_t^X$ forms a lattice under the partial order given by \eqref{eqn_partial_order}.
\end{proposition}

\begin{proof}
Closure under $\min$ and $\max$ follows from the Proposition \ref{proposition_ptwise_algebraic_structure}. The partial order is the standard pointwise order on random variables given by \eqref{eqn_partial_order}. The lattice properties follow from the fact that $\min$ and $\max$ define meet and join operations described by \eqref{eqn_lattices} under this order.
\end{proof}

\subsubsection{Order Preservation}
\noindent We next relate logical implication to the ordering of observables.
\begin{proposition}
Let $\phi_1, \phi_2 \in \Phi$ be such that $\rho_{\phi_1}(x,t) \leq \rho_{\phi_2}(x,t), \forall x \in \mathcal X^{\mathbb{T}}$. Then, combining the formulae with the stochastic process $X$, we see that $S_{\phi_1}^X(t) \leq S_{\phi_2}^X(t)$
almost surely.
\end{proposition}

\begin{proof}
For every $\omega\in\Omega$, we observe that $S_{\phi_1}^X(t)(\omega) = \rho_{\phi_1}(X(\omega),t) \leq \rho_{\phi_2}(X(\omega),t) = S_{\phi_2}^X(t)(\omega)$. Thus the inequality holds pointwise and hence almost surely.
\end{proof}

\textit{Example:} Consider the two temporal formulas defined on the stochastic process $X$ namely $\phi_1 := (X_t \geq 10)$, and $\phi_5 := (X_t \geq 5)$. The predicate $\phi_1 := (X_t \geq 10)$ is stricter than $\phi_1 := (X_t \geq 5)$. Their corresponding robustness values are $\rho_{\phi_{1}}(x, t) := x_t - 10$, and $\rho_{\phi_{2}}(x, t) := x_t - 5$ respectively. Clearly, $x_t - 10 \leq x_t - 5$ meaning that satisfying $X_t \geq 10$ is harder than satisfying $X_t \geq 5$. Hence, $S_{\phi_1}^X(t) \leq S_{\phi_2}^X(t)$ meaning that the stochastic observable preserves this ordering.

\begin{remark}
The above exposition establishes a bridge between temporal logic and lattice-valued stochastic analysis, enabling the use of algebraic and order-theoretic tools in the study of temporal properties of stochastic processes.
\end{remark}
\section{Statistical Analysis of Temporal Satisfaction} \label{sec_5}
In this section, we study how can ordinary statistics be applied to study temporal satisfaction of a formula by a stochastic process $X$. For a finite $N \in \mathbb{N}$, let $\{\omega_i\}_{i=1}^N \subset \Omega$ denote the independent samples drawn according to distribution $\mathbb{P}$. Then, the empirical satisfaction average at time $t$ is defined as
\begin{align}
\label{eqn_empirical_satisfaction_average}
\widehat{p}_N(t)
:=
\frac{1}{N}\sum_{i=1}^N
\mathds{1}_{\phi}^X(t)(\omega_i).
\end{align}
\begin{definition}
For $\epsilon > 0$, the $\epsilon$-deviation event at time $t$ is defined as
\begin{align}
\label{eqn_deviation_event}
A_{N,\epsilon}(t)
:=
\left\{
\left\{ \omega_i \right\}^{N}_{i=1} 
:
\left|
\widehat{p}_N(t)
-
\mathbb{P}\big(X,t \models\phi\big)
\right|
\geq 
\epsilon
\right\}.
\end{align}
\end{definition}
\noindent If $\mathbb{P}\big(S_\phi^X(t)=0\big) = 0$, then $\mathbb{P}\big(X,t \models\phi\big) = \mathbb{P}\big(S_\phi^X(t)\geq0\big)$ and hence the empirical satisfaction average also estimates the
probability that the TLOP is nonnegative. We will now prove the law of large numbers for temporal satisfaction in the following lemma.
\begin{lemma}\label{lemma_lln_temporal_satisfy}
$\widehat p_N(t) \xrightarrow[N\to\infty]{\mathrm{a.s.}}
\mathbb{P}(X,t \models \phi)$.
\end{lemma}

\begin{proof}
Define $Y_i(t) := \mathds{1}_\phi^X(t)(\omega_i)$, for $i \in \mathbb{N}$. Since $\omega_i$ are independent samples drawn according to $\mathbb{P}$, the
random variables $\{Y_i(t)\}_{i \in \mathbb{N}}$ are independent and identically
distributed Bernoulli random variables. Moreover,
\begin{align*}
\mathbb{E}[Y_i(t)] = \mathbb{P}(X,t \models \phi).
\end{align*}
Since $Y_i(t) \in [0,1]$, we have $\mathbb{E}[|Y_i(t)|]<\infty$. Hence,
by the strong law of large numbers, the result follows.
\end{proof}

\begin{corollary}
The Hoeffding concentration for temporal satisfaction $\forall \epsilon>0$ is given by
\begin{align}
\label{eqn_corollary_Hoeffding}
\mathbb{P}\bigl(A_{N,\epsilon}(t)\bigr)
\leq
2\exp(-2N\epsilon^2).
\end{align}
\end{corollary}

\begin{proof}
From Lemma \ref{lemma_lln_temporal_satisfy}, we know that $\{Y_i(t)\}_{i=1}^N$ are independent random variables with $Y_i(t) := \mathds{1}_\phi^X(t)(\omega_i)$ satisfying $Y_i(t) \in [0,1]$.
By Hoeffding's inequality for independent bounded random variables in \citep{vershynin2018high}, we see that
\begin{align*}
\mathbb{P}\left(
\left|
\widehat p_N(t)
-
\mathbb{E}[\widehat p_N(t)]
\right|
\geq \epsilon
\right)
\leq
2\exp(-2N\epsilon^2).
\end{align*}
The result \eqref{eqn_corollary_Hoeffding} follows immediately by using the definition \eqref{eqn_deviation_event} and the fact that $\mathbb{E}[\widehat p_N(t)] = \mathbb{P}(X,t \models \phi)$.
\end{proof}

Consider the same stock market example explained before. If independent simulated market trajectories
$\omega_1,\ldots,\omega_N$ are generated, for example by simulating independent
Brownian paths, the empirical satisfaction average $\widehat p_N(t)$ is given by \eqref{eqn_empirical_satisfaction_average}. This quantity estimates the probability that the stock satisfies the temporal
formula. The law of large numbers then provides a probabilistic justification
for this empirical estimate, while concentration inequalities quantify the
finite-sample deviation between the empirical estimate and the true
satisfaction probability.


\section{Regularity of Temporal-Logical Observables} \label{sec_6}

In this section, we study the pathwise regularity of temporal-logical
observable processes. The measurability results established earlier show that,
for each fixed time $t \in \mathbb{T}$, the temporal-logical observable $S_\phi^X(t)$ is a well-defined random variable. We now ask whether the entire process $\{S_\phi^X(t)\}_{t\geq 0}$ inherits regularity properties from the underlying stochastic process $X=\{X_t\}_{t\geq 0}$. We focus on the continuous-time case $\mathbb{T}=\mathbb R_{\geq 0}$. The main result of this section shows that, if the sample paths of $X$ are continuous and the atomic predicates are defined by continuous functions, then the TLOP has continuous sample paths, provided all temporal intervals are compact. For the purposes of this paper, we restrict attention to finite-valued robustness semantics so that all temporal-logical observables take values in $\mathbb R$. This avoids technical
complications arising from extended real-valued random variables. To this end, we use finite robustness values for the Boolean constants. That is, for a fixed constant $M>0$, we set $\rho_{\top}(x,t) := M$, and $\rho_{\bot}(x,t) := - M$.

\begin{theorem}\label{thm_pathwise_continuity}
Let $\mathbb{T}=\mathbb R_{\geq 0}$, $I=[a,b]$ with $0\leq a\leq b<\infty$, and $\mathcal X\subseteq \mathbb R^d$. Let $X = \{X_t\}_{t\geq 0}$ be an $\mathcal X$-valued stochastic process on $(\Omega,\mathcal F,\mathbb{P})$ with continuous sample paths almost surely. Further, let $h_\mu:\mathcal X\to\mathbb R$ defining every atomic predicate $\mu$ be a continuous function. Then, for almost every
$\omega\in\Omega$, the mapping $t\mapsto S_\phi^X(t)(\omega)$ is continuous on $\mathbb R_{\geq 0}$.
\end{theorem}

\begin{proof}
Let $\Omega_c\subseteq \Omega$ denote the event on which the sample path
$t\mapsto X_t(\omega)$ is continuous. By assumption, $\mathbb{P}(\Omega_c) = 1$. Fix $\omega\in\Omega_c$. We will now prove that, for this fixed sample path
$X(\omega)$, the mapping $t \mapsto S_\phi^X(t)(\omega) = \rho_\phi(X(\omega),t)$ is continuous. The result then holds almost surely. As before, the proof will proceed by structural induction on the formula $\phi$. \\

\noindent Case 1: Let $\phi=\mu$. Then $S_\mu^X(t)(\omega) =
h_\mu(X_t(\omega))$. Since $t\mapsto X_t(\omega)$ is continuous and $h_\mu$ is continuous, the composition $t \mapsto h_\mu(X_t(\omega))$ is continuous. Hence $t\mapsto S_\mu^X(t)(\omega)$ is continuous.

\noindent Case 2: Assume that the result holds for a formula $\psi$. Let $\phi=\neg\psi$. By the quantitative semantics, $S_{\neg\psi}^X(t)(\omega) =
- S_\psi^X(t)(\omega)$. By the induction hypothesis, $t\mapsto S_\psi^X(t)(\omega)$ is continuous. Therefore its negative is continuous, and hence
$t\mapsto S_{\neg\psi}^X(t)(\omega)$ is continuous. 

\noindent Case 3: Assume that the result holds for formulas $\psi_1$ and $\psi_2$.
For conjunction and disjunction
\begin{align*}
S_{\psi_1\wedge\psi_2}^X(t)(\omega)
&=
\min\{S_{\psi_1}^X(t)(\omega),S_{\psi_2}^X(t)(\omega)\}, \quad \text{and} \\
S_{\psi_1\vee\psi_2}^X(t)(\omega)
&=
\max\{S_{\psi_1}^X(t)(\omega),S_{\psi_2}^X(t)(\omega)\}
\end{align*}
respectively. Since the pointwise minimum of two continuous real-valued functions is continuous and the pointwise maximum of two continuous real-valued functions is continuous, we observe that both $t \mapsto S_{\psi_1\wedge\psi_2}^X(t)(\omega)$ and $t \mapsto S_{\psi_1\vee\psi_2}^X(t)(\omega)$ are continuous.

\noindent Case 4: Assume that the result holds for a formula $\psi$. Let
$\phi=\Box_{I}\psi$, where $I=[a,b]$ is compact. Define $r(t):=S_\psi^X(t)(\omega)$. By the induction hypothesis, $r:\mathbb R_{\geq 0}\to\mathbb R$ is continuous. The robustness of the always formula is
\begin{align*}
S_{\Box_{I}\psi}^X(t)(\omega)
=
\inf_{s\in[t+a,t+b]} r(s).
\end{align*}
We claim that the mapping $g(t) := \inf_{s\in[t+a,t+b]} r(s)$ is continuous. Fix $t_0\geq 0$ and choose $\delta>0$. For all $t$ sufficiently close to
$t_0$, the intervals $[t+a,t+b]$ are contained in a common compact interval
$K_{t_0} := [\max\{0,t_0+a-\delta\},t_0+b+\delta]$. Since $r$ is continuous, it is uniformly continuous on $K_{t_0}$. Let $\epsilon > 0$. By uniform continuity, there exists $\eta>0$ such that whenever $u,v\in K_{t_0}$ and $|u-v|<\eta$, we have $|r(u)-r(v)|<\epsilon$. Now take $t$ such that $|t-t_0|<\eta$ and $|t-t_0|<\delta$. For every $s \in [t+a,t+b]$, define $s_0 := s - (t - t_0)$. Then $s_0 \in [t_0 + a, t_0 + b]$ and $|s - s_0| = |t - t_0| < \eta$. Hence, $r(s)\geq r(s_0)-\epsilon$. Taking the infimum over $s\in[t+a,t+b]$ gives $g(t) \geq g(t_0) - \epsilon$. By reversing the roles of $t$ and $t_0$, we similarly obtain $g(t_0) \geq g(t) - \epsilon$. Therefore, $|g(t)-g(t_0)|\leq \epsilon$. Since $\epsilon>0$ was arbitrary, $g$ is continuous at $t_0$. As $t_0$ was arbitrary, $g$ is continuous on $\mathbb R_{\geq 0}$. Thus $t \mapsto S_{\Box_{I}\psi}^X(t)(\omega)$ is continuous.

\noindent Case 5: The proof is analogous for the eventually operator. Let $\phi=\Diamond_{I}\psi$, where $I=[a,b]$ is compact, and again define $r(t) := S_\psi^X(t)(\omega)$. By the induction hypothesis, $r$ is continuous. The robustness of the eventually formula is
\begin{align*}
S_{\Diamond_{I}\psi}^X(t)(\omega)
=
\sup_{s\in[t+a,t+b]} r(s).
\end{align*}
Define $q(t):=\sup_{s\in[t+a,t+b]} r(s)$. Using the same uniform-continuity argument as in the previous step, one obtains $|q(t)-q(t_0)|\leq \epsilon$ whenever $t$ is sufficiently close to $t_0$. Hence $q$ is continuous. Therefore $t\mapsto S_{\Diamond_{I}\psi}^X(t)(\omega)$ is continuous. \\

\noindent Since every formula is obtained by finitely many applications of the grammar rules, the result follows by structural induction. Therefore, for every
$\omega\in\Omega_c$, the path $t\mapsto S_\phi^X(t)(\omega)$ is continuous. Since $\mathbb{P}(\Omega_c)=1$, the TLOP $S_\phi^X$ has continuous sample paths almost surely.
\end{proof}

\section{Distributional Properties and Extremal Characterization} \label{sec_7}

In this section, we investigate the distributional properties of
temporal-logical observables. In particular, we show that, for a broad class
of temporal formulas, the distribution of $S_\phi^X(t)$ can be characterized
in terms of extremal functionals of the underlying stochastic process $X$. The key observation is that temporal operators such as always and
eventually naturally correspond to infimum and supremum operations over time intervals. As a result, temporal-logical observables
encode information about extrema of the process, linking TLSP directly to
classical problems in stochastic analysis such as barrier crossing and
hitting probabilities.

\subsection{Extremal Characterization of Temporal Satisfaction}
Given a fixed constant $c\in\mathbb R$, we focus on temporal formulas constructed from threshold predicates of the
form $\mu(x) := x - c$.

\begin{theorem}
Let $X=\{X_t\}_{t\geq 0}$ be a real-valued stochastic process with continuous sample paths, and let $t\geq 0$ and $T>0$ be fixed. Let $\phi_1 = \Box_{[0,T]}(X \geq c)$, $\phi_2 = \Diamond_{[0,T]}(X \geq c)$, and $\phi_3 = \Diamond_{[0,T]}(X \geq K) \wedge \Box_{[0,T]}(X \geq L)$, with $K\geq L$. Then,
\begin{align}
\mathbb{P}(X,t \models \phi_1)
&=
\mathbb{P}\Big(\inf_{s\in[t,t+T]} X_s \geq c\Big), \label{eqn_thm_10_result_1}\\
\mathbb{P}(X,t \models \phi_2)
&=
\mathbb{P}\Big(\sup_{s\in[t,t+T]} X_s \geq c\Big), \label{eqn_thm_10_result_2}\\
\mathbb{P}(X,t \models \phi_3)
&=
\mathbb{P}\Big(
\sup_{s\in[t,t+T]} X_s \geq K
\;\text{and}\;
\inf_{s\in[t,t+T]} X_s \geq L
\Big). \label{eqn_thm_10_result_3}
\end{align}
\end{theorem}

\begin{proof}
We prove each statement directly from the Boolean semantics of the
temporal operators. Since $X$ has continuous sample paths and $[t,t+T]$ is compact, the minimum and maximum of $X_s$ over this interval are attained almost surely. Given $\phi_1 = \Box_{[0,T]}(X \geq c)$, we see from \eqref{eqn_temporal_logical_observable} that, 
\begin{align*}
X,t \models \phi_1
\iff
X_s \geq c, \forall s\in[t,t+T]
\iff 
\inf_{s\in[t,t+T]} X_s \geq c.
\end{align*}
Taking probabilities on both sides yields \eqref{eqn_thm_10_result_1}. Now consider $\phi_2 = \Diamond_{[0,T]}(X \geq c)$. Then, by similar argument
\begin{align*}
X,t \models \phi_2
\iff 
\exists s\in[t,t+T] \text{ such that }
X_s \geq c
\iff
\sup_{s\in[t,t+T]} X_s \geq c.
\end{align*}
Taking probabilities on both sides yields \eqref{eqn_thm_10_result_2}. Let
$\phi_3 = \Diamond_{[0,T]}(X \geq K) \;\wedge\; \Box_{[0,T]}(X \geq L)$. 
By the Boolean semantics of conjunction, eventually, and always,
\begin{align*}
(X,t)\models\phi_3
&\Longleftrightarrow
\left(
\exists s\in[t,t+T]\text{ such that }X_s\geq K
\right)
\wedge
\left(
X_s\geq L,\ \forall s\in[t,t+T]
\right) \\
&\Longleftrightarrow
\left\{
\sup_{s\in[t,t+T]}X_s\geq K
\right\}
\cap
\left\{
\inf_{s\in[t,t+T]}X_s\geq L
\right\}.
\end{align*}
Finally, taking probability on both sides, yields \eqref{eqn_thm_10_result_3}.
\end{proof}
Having learned the way to turn temporal logic problems into extremal stochastic-process problems, we will now use the above developments to study about explicit TLSP satisfaction probabilities for Brownian motion.
\begin{theorem}
\label{thm_reflection_principle}
Let $\Phi_{\mathcal{N}}$ denote the standard normal distribution function and further let $W=\{W_s\}_{s\geq 0}$ be a standard Brownian motion and define, for
$s\geq t$,
\begin{align} \label{eqn_Brownian_motion}
X_s = x+\sigma(W_s-W_t),
\end{align}
where $x\in\mathbb R$ and $\sigma>0$. Let $T>0$. Define the eventual-reaching formula and the safety formula respectively as $\phi_K := \Diamond_{[0,T]}(X\geq K)$, and $\psi_L := \Box_{[0,T]}(X\geq L)$. Then, 
\begin{align}
\mathbb{P}(X,t \models \phi_K\mid X_t=x)
&=
\begin{cases}
1, & x\geq K,\\
2\left(1-\Phi_{\mathcal{N}}\left(\dfrac{K-x}{\sigma\sqrt T}\right)\right),
& x<K,
\end{cases}, \quad \text{and} \label{eqn_eventual_reach_satisfaction_probability}\\
\mathbb{P}(X,t \models \psi_L\mid X_t=x)
&=
\begin{cases}
0, & x<L,\\
2\Phi_{\mathcal{N}}\left(\dfrac{x-L}{\sigma\sqrt T}\right)-1,
& x\geq L. \label{eqn_always_safe_satisfaction_probability}
\end{cases}
\end{align}
\end{theorem}

\begin{proof}
We first prove the formula for $\phi_K$. By the definition of $\Diamond_{[0,T]}$, $(X,t) \models \phi_K \iff \sup_{s\in[t,t+T]} X_s \geq K$. If $x\geq K$, the event holds immediately, and the probability is $1$. Assume now that $x<K$. From \eqref{eqn_Brownian_motion}, we have
\begin{align}
\sup_{s\in[t,t+T]}X_s\geq K
\iff
\sup_{0\le u\le T}(W_{t+u}-W_t)
\geq
\frac{K-x}{\sigma}.
\end{align}
By stationary independent increments of Brownian motion, the process, $\widetilde W_u:=W_{t+u}-W_t, 0\le u\le T$, is again a standard Brownian motion. Hence
\begin{align}
\mathbb{P}(X,t \models \phi_K\mid X_t=x)
=
\mathbb{P}\left(
\sup_{0\le u\le T}\widetilde W_u
\geq
\frac{K-x}{\sigma}
\right).
\end{align}
By the reflection principle,
\begin{align}
\mathbb{P}\left(
\sup_{0\le u\le T}\widetilde W_u\geq a
\right)
=
2\left(1-\Phi_{\mathcal{N}}\left(\frac{a}{\sqrt T}\right)\right),
\qquad a>0.
\end{align}
Taking $a=(K-x)/\sigma$ gives \eqref{eqn_eventual_reach_satisfaction_probability}. We now prove the formula for $\psi_L$. By the definition of
$\Box_{[0,T]}$, $(X,t) \models \psi_L \iff \inf_{s\in[t,t+T]} X_s \geq L$. If $x<L$, the event is impossible at the initial time, so the probability is $0$. Assume now that $x\geq L$ and let $a:=\frac{x-L}{\sigma}\geq 0$. Then
\begin{align}
\inf_{s\in[t,t+T]}X_s\geq L
\iff
\inf_{0\le u\le T}\widetilde W_u
\geq
-a
\iff
\sup_{0\le u\le T}(-\widetilde W_u)\le a.
\end{align}
Since $-\widetilde W$ is again a standard Brownian motion, the reflection
principle gives
\begin{align}
\mathbb{P}\left(
\sup_{0\le u\le T}(-\widetilde W_u)\geq a
\right)
=
2\left(1-\Phi_{\mathcal{N}}\left(\frac{a}{\sqrt T}\right)\right).
\end{align}
Therefore,
\begin{align}
\mathbb{P}\left(
\sup_{0\le u\le T}(-\widetilde W_u)< a
\right)
=
1-
2\left(1-\Phi_{\mathcal{N}}\left(\frac{a}{\sqrt T}\right)\right).
\end{align}
For Brownian motion, the distinction between $<$ and $\leq$ at the fixed barrier level does not affect the probability. Hence
the result \eqref{eqn_always_safe_satisfaction_probability} follows and this completes the proof.
\end{proof}

\subsection{First Hitting-Time Characterization of TLSP Satisfaction}
\begin{definition}
Given a real-valued stochastic process $X=\{X_s\}_{s\geq t}$ with continuous sample paths, the first upper hitting time and the first lower violation time with respective to levels $K, L \in \mathbb{R}$ are respectively defined as
\begin{align}
\tau_K^+
&:=
\inf\{s \geq 0 \mid X_{t+s} \geq K\}, \quad \text{and}\\
\tau_L^-
&:=
\inf\{s \geq 0 \mid X_{t+s} < L\}.
\end{align}    
\end{definition}

\begin{theorem}\label{thm_hitting_time}
Fix $t \geq 0$. Let $T > 0$, and a real-valued stochastic process $X = \{X_s\}_{s\ge t}$ with continuous sample paths. Then, almost surely,
\begin{align}
    X, t \models \lozenge_{[0,T]}(X\ge K)
    &\Longleftrightarrow
    \tau_K^+\le T,\\
    X, t \models \square_{[0,T]}(X\ge L)
    &\Longleftrightarrow
    \tau_L^-\ge T.
\end{align}
Consequently, for the combined formula, we have
\begin{align}
    \label{eqn_thm_hitting_time_result}
    \mathbb{P}
    \left(
    X, t 
    \,
    \models 
    \,
    \underbrace{\lozenge_{[0,T]}(X\ge K) \wedge \square_{[0,T]}(X\ge L)}_{:= \phi}
    \right)
    =
    \mathbb{P}
    \left(
        \tau_K^+\le T,\,
        \tau_L^-\ge T
    \right).
\end{align}
\end{theorem}
\begin{proof}
Let $\Omega_c$ denote the probability-one event on which the sample paths of $X$ are continuous. Fix $\omega \in \Omega_c$. We first consider the eventually formula $\varphi_K = \Diamond_{[0,T]}(X\geq K)$. By the Boolean semantics of the eventually operator,
\begin{align*}
(X(\omega),t)\models\varphi_K 
\iff \exists u\in[0,T] \mid X_{t+u}(\omega)\geq K.
\end{align*}
Since the sample path is continuous, the set $\{u\geq0 \mid X_{t+u}(\omega)\geq K\}$ is closed. Whenever this set is non-empty, its infimum is attained. Therefore,
\begin{equation*}
(X(\omega),t)\models\varphi_K
\iff \tau_K^+(\omega)\leq T.
\end{equation*}
Next, consider the always formula $\psi_L=\Box_{[0,T]}(X\geq L)$. By the Boolean semantics of the always
operator, 
\begin{align*}
(X(\omega),t)\models\psi_L 
\iff X_{t+u}(\omega)\geq L,
\quad\forall u\in[0,T].
\end{align*}
If $\psi_L$ is satisfied, there is no strict violation $X_{t+u}(\omega)<L$ for
$u\in[0,T]$. Hence, $\tau_L^-(\omega)\geq T$. Conversely, suppose that $\tau_L^-(\omega)\geq T$, but that
$\psi_L$ is violated. Then there exists $u_0\in[0,T]$ such that $X_{t+u_0}(\omega)<L$. If $u_0<T$, then immediately $\tau_L^-(\omega)\leq u_0<T$, contradicting the assumption. If $u_0=T$, continuity of the sample path implies that $\exists \delta>0$ such that $X_{t+u}(\omega)<L$, for 
$u \in (T-\delta,T]$, where $\delta$ may be chosen smaller than $T$. Consequently, a strict violation occurs
before $T$, implying $\tau_L^-(\omega)<T$, which is again a contradiction.
Hence,
\begin{equation*}
(X(\omega),t)\models\psi_L
\iff \tau_L^-(\omega)\geq T.
\end{equation*}
Finally, consider the combined formula $\varphi=\varphi_K\wedge\psi_L$. By the Boolean semantics of conjunction,
\begin{align*}
(X(\omega),t)\models\varphi
\iff
\left(\tau_K^+(\omega)\leq T\right)
\wedge
\left(\tau_L^-(\omega)\geq T\right).
\end{align*}
Since this equivalence holds $\forall \omega\in\Omega_c$ and
$\mathbb{P}(\Omega_c)=1$, taking probabilities yields the result \eqref{eqn_thm_hitting_time_result} and the proof is complete. 
\end{proof}
Note that the above exposition nicely connects with the Doob h-transforms and the corresponding investigation is left as a future work.

\section{Perturbation and Robustness of Temporal Satisfaction} \label{sec_8}
In this section, we study the sensitivity of temporal-logical observables with respect to perturbations in the underlying stochastic process. This is important in practical settings where the process is subject to model uncertainty, approximation error, or measurement noise. Our goal is to establish that temporal satisfaction defined using temporal robustness is stable under uniform perturbations of sample paths. In particular, we show that temporal-logical observables are Lipschitz continuous with respect to the supremum norm on trajectories. To this end, we consider only real-valued processes defined on $\mathbb{T}=\mathbb R_{\geq 0}$ meaning that $\mathcal X\subseteq \mathbb R$.

\begin{assumption} \label{assume_h_mu_lipschitz}
For each atomic predicate $\mu$, $\exists L_\mu>0$ such that the Borel measurable function $h_\mu:\mathcal X \to \mathbb R$ defining it is $L_\mu$-Lipschitz continuous. That is,
\begin{align}
|h_\mu(x)-h_\mu(y)| \le L_\mu |x-y|,
\quad \forall x,y\in\mathcal X.
\end{align}
\end{assumption}

\begin{lemma} \label{lemma_real_numbers_relation}
For any $a,b,c,d \in \mathbb{R}$, the following results hold:
\begin{align}
|\min(a,b)-\min(c,d)|
&\leq
\max\{|a-c|,|b-d|\}, \label{eqn_lemma_real_numbers_relation_1} \\
|\max(a,b)-\max(c,d)|
&\leq
\max\{|a-c|,|b-d|\}.
\label{eqn_lemma_real_numbers_relation_2}
\end{align}    
\end{lemma}
\begin{lemma} \label{lemma_functions_relation}
For any two functions $f,g$, the following results hold:
\begin{align}
|\inf f - \inf g|
&\leq
\sup |f-g|, \label{eqn_lemma_functions_relation_1}\\
|\sup f - \sup g|
&\leq
\sup |f-g|. \label{eqn_lemma_functions_relation_2}
\end{align}
\end{lemma}

\subsection{Perturbation Studies of Temporal-Logical Observables}
Robustness of temporal logic specifications for continuous time signals have been previously studied in \citep{fainekos2009robustness}, albeit in the deterministic setting. However, we elevate such an exposition to the TLOP setting with temporal robustness becoming a random variable for every time in the following theorem.
\begin{theorem} \label{thm_pathwise_robustness}
Given two processes $X=\{X_t\}_{t\geq 0}$ and $\tilde X=\{\tilde X_t\}_{t\geq 0}$ defined on $(\Omega,\mathcal{F},\mathbb{P})$ with continuous sample paths, let each interval $I$ defining the temporal formula $\phi \in \Phi$ be compact. Then, under assumption \ref{assume_h_mu_lipschitz}, $\exists C_\phi>0$, depending only on the structure of $\phi$, $L_\mu$, and an event $\Omega_0\in\mathcal{F}$ with
$\mathbb{P}(\Omega_0)=1$, such that $\forall t \geq 0$ and every sample path $\omega \in \Omega_0$,
\begin{align}
\label{eqn_thm15_result}
|S_\phi^X(t)(\omega) - S_\phi^{\tilde X}(t)(\omega)|
\le
C_\phi
\sup_{s\in[t,t+H_\phi]}
|X_s(\omega)-\tilde X_s(\omega)|.
\end{align}
\end{theorem}

\begin{proof}
Let $\Omega_X$ and $\Omega_{\widetilde X}$ denote the certainty events on which the sample paths of $X$ and $\widetilde X$, respectively, are continuous. Define $\Omega_0:=\Omega_X\cap\Omega_{\widetilde X}$. Then $\mathbb{P}(\Omega_0)=1$. Fix $\omega\in\Omega_0$. The proof then  proceeds by structural induction on $\phi$.

\noindent Case 1: Consider the setting where the formula is given by atomic predicates meaning that 
$\phi=\mu$. Then, 
$S_\mu^X(t)(\omega) =
h_\mu(X_t(\omega))$, and
$S_\mu^{\tilde X}(t)(\omega)
=
h_\mu(\tilde X_t(\omega))$. By Lipschitz continuity,
\begin{align*}
|S_\mu^X(t)-S_\mu^{\tilde X}(t)|
\le
L_\mu |X_t-\tilde X_t|
\le
L_\mu \sup_{s\in[t,t+H_\phi]}|X_s-\tilde X_s|.
\end{align*}

\noindent Case 2: Consider the setting where the formula is given by negation meaning that $\phi = \neg\psi$. Then $S_{\neg\psi}^X(t) = - S_\psi^X(t)$. Thus, $|S_{\neg\psi}^X(t) -S_{\neg\psi}^{\tilde X}(t)| =
|S_\psi^X(t) - S_\psi^{\tilde X}(t)|$. The bound follows from the induction hypothesis.

\noindent Case 3: Consider the setting where the formula is given by either conjunction or disjunction. Let $\phi=\psi_1\wedge\psi_2$. Then, $S_\phi^X(t) = \min \{ S_{\psi_1}^X(t), S_{\psi_2}^X(t)\}$. 
Applying \eqref{eqn_lemma_real_numbers_relation_1}, we get
\begin{align*}
|S_\phi^X(t)-S_\phi^{\tilde X}(t)|
\leq
\max\{
|S_{\psi_1}^X(t)-S_{\psi_1}^{\tilde X}(t)|,
|S_{\psi_2}^X(t)-S_{\psi_2}^{\tilde X}(t)|
\}.
\end{align*}
Using the induction hypothesis for $\psi_1$ and $\psi_2$ gives the result. The case $\phi=\psi_1\vee\psi_2$ follows similarly by applying \eqref{eqn_lemma_real_numbers_relation_2}. 

\noindent Case 3: Consider the setting where the formula is given by always operator meaning that $\phi = \Box_{I}\psi$. Then $S_\phi^X(t) = \inf_{s\in[t,t+I]} S_\psi^X(s)$. Using \eqref{eqn_lemma_functions_relation_1}, we get
\begin{align*}
|S_\phi^X(t)-S_\phi^{\tilde X}(t)|
\le
\sup_{s\in[t,t+I]}
|S_\psi^X(s)-S_\psi^{\tilde X}(s)|.
\end{align*}
Applying the induction hypothesis yields the bound.

\noindent Case 4: Consider the setting where the formula is given by the eventually operator meaning that $\phi = \Diamond_{I}\psi$.
The proof is identical to the always operator, using \eqref{eqn_lemma_functions_relation_2}. Combining all steps, we get a constant $C_\phi$ depending only on the formula structure and predicate Lipschitz constants, such that \eqref{eqn_thm15_result} holds true. That is, $C_\mu = L_\mu, C_{\neg\psi} = C_\psi, C_{\psi_1\wedge\psi_2}
= C_{\psi_1\vee\psi_2}
= \max\{C_{\psi_1},C_{\psi_2}\}$, and $C_{\square_I\psi} = C_{\lozenge_I\psi} = C_\psi$. This completes the proof.
\end{proof}
We just proved that the robustness transformation $X \mapsto S_{\phi}^{X}$ is Lipschitz with respect to the supremum norm. Informally, if the two processes satisfy $X \approx \tilde{X}$, then their associated temporal robustness satisfy $S_{\phi}^{X} \approx S_{\phi}^{\tilde{X}}$.

\subsection{Probability-Level Robustness}
The previous result gives a pathwise bound between temporal-logical observables. Now consider that at time $t$, the stochastic process $X$ satisfying the formula $\phi$ results in $S^X_{\phi}(t) = 100$. Then, a smaller perturbation, say $0.1$ will not affect the sign of $S^X_{\phi}(t)$ and hence the probability of satisfaction of formula $\phi$ by the underlying process $X$. On the contrary, now suppose that $S^X_{\phi}(t) = 10^{-5}$. Then,  same perturbation of $0.1$ can now flip the formula satisfaction to violation. So, pathwise robustness does not automatically imply satisfaction-probability robustness. Probability-level robustness depends not only on how close $X$ and $\widetilde X$ are, but also on how much probability mass lies near $S^X_{\phi}(t) = 0$. \\
We must derive a corresponding probability-level robustness result. The key point is that satisfaction probabilities are stable under perturbations provided that the original temporal-logical observable has small probability
mass near the satisfaction boundary. We now study probability-level robustness of temporal satisfaction in the following theorem which says that there are two ways through which formula satisfaction can change namely when, 1) the process perturbation is too large, and 2) the original trajectory of the process was already close to the logical boundary of the temporal satisfaction.

\begin{theorem}\label{thm_probability_satisfaction_robustness}
Let $X=\{X_s\}_{s\geq 0}$ and $\tilde X=\{\tilde X_s\}_{s\geq 0}$ be two real-valued stochastic processes defined on the common probability space \((\Omega,\mathcal F,\mathbb P)\), satisfying the assumptions of Theorem \ref{thm_pathwise_robustness}. Let $\phi$ satisfy Theorem \ref{thm_pathwise_robustness} with constant $C_\phi>0$. For fixed $t\geq 0$, define the pathwise deviation over $H_\phi$ of $\phi$ by
\begin{align}
\label{eqn_thm16_assume_1}
\Delta_t
:=
\sup_{s\in[t,t+H_\phi]}
|X_s-\tilde X_s|.
\end{align}
Then, $\forall \eta>0$,
\begin{align} \label{eqn_thm16_result}
\left|
\mathbb{P}(X,t \models\phi)
-
\mathbb{P}( (\tilde X,t) \models\phi)
\right|
\leq
\mathbb{P}\left(
|S_\phi^X(t)|\leq C_\phi\eta
\right)
+
\mathbb{P}(\Delta_t>\eta).
\end{align}
\end{theorem}

\begin{proof}
For a fixed $t\geq 0$, define the Boolean satisfaction events $A := \{\omega\in\Omega \mid(X(\omega),t)\models\varphi\}$, and  
$\widetilde A := \{\omega\in\Omega \mid
(\widetilde X(\omega),t)\models\varphi\}$.
By definition, $\mathbb{P}(X,t\models\varphi)
= \mathbb{P}(A)$, and 
$\mathbb{P}(\widetilde X,t\models\varphi) =
\mathbb{P}(\widetilde A)$. Let $A \triangle \widetilde A$ denote the symmetric
difference between the two events. Then,
\begin{equation}
\left|
\mathbb{P}(X,t\models\varphi)
-
\mathbb{P}(\widetilde X,t\models\varphi)
\right|
\leq
\mathbb{P}(A\triangle\widetilde A).
\end{equation}
We now establish that, up to a null event,
\begin{equation}
A\triangle\widetilde A
\subseteq
\left\{
|S_\varphi^X(t)|\leq C_\varphi\eta
\right\}
\cup
\left\{\Delta_t>\eta\right\}.
\end{equation}
To this end, consider an outcome $\omega$ belonging
to the probability-one event on which the pathwise
bound in Theorem \ref{thm_pathwise_robustness} holds. Suppose that $\Delta_t(\omega) \leq \eta$. By Theorem \ref{thm_pathwise_robustness}, we have
\begin{equation}
\left|
S_\varphi^X(t)(\omega)
-
S_\varphi^{\widetilde X}(t)(\omega)
\right|
\leq C_\varphi\eta.
\end{equation}
Now suppose, in addition, that
$|S_\varphi^X(t)(\omega)|>C_\varphi\eta$.
There are two possible cases.
\begin{enumerate}
    \item \textit{Case 1:} If $S_\varphi^X(t) (\omega) > C_\varphi\eta$, then
    \begin{align*}
    S_\varphi^{\widetilde X}(t)(\omega)
    &\geq
    S_\varphi^X(t)(\omega)-C_\varphi\eta
    >
    0.
    \end{align*}
    Thus, both temporal-logical observables have strictly positive robustness. By the sign-soundness property of quantitative semantics, both $X(\omega)$ and $\widetilde X(\omega)$ satisfy $\varphi$ at time $t$. Consequently, $\omega \in A\cap\widetilde A$.
    \item \textit{Case 2:} If $S_\varphi^X(t)(\omega)<-C_\varphi\eta$, then
    \begin{align*}
    S_\varphi^{\widetilde X}(t)(\omega)
    &\leq
    S_\varphi^X(t)(\omega)+C_\varphi\eta
    <0.
    \end{align*}
    Thus, both temporal-logical observables have strictly negative robustness. By sign soundness, both $X(\omega)$ and $\widetilde X(\omega)$ violate $\varphi$ at time $t$. Consequently, $\omega \notin A\cup\widetilde A$. 
\end{enumerate}
In either case, the Boolean satisfaction status of the two processes agrees, implying
$\omega\notin A\triangle\widetilde A$. Therefore, outside a null event, if $\omega\in A\triangle\widetilde A$, then either
$|S_\varphi^X(t)(\omega)|\leq C_\varphi\eta$ or $\Delta_t(\omega)>\eta$. Hence,
\begin{equation}
\mathbb{P}(A\triangle\widetilde A)
\leq
\mathbb{P}\left(
|S_\varphi^X(t)|\leq C_\varphi\eta
\right)
+
\mathbb{P}(\Delta_t>\eta).
\end{equation}
Combining this bound with the preceding
symmetric-difference inequality yields \eqref{eqn_thm16_result} and the proof is complete. 
\end{proof}
The following corollary establishes the robustness under almost sure uniform perturbations. That is, if the process perturbation is negligible, then the satisfaction can change only when the original trajectory of the process was already close to the logical boundary of the temporal satisfaction. 
\begin{corollary}
Under the assumptions of the previous theorem, suppose that $\Delta_t \leq \eta$ almost surely.  Then
\begin{align}
\label{corollary_16_result}
\left|
\mathbb{P}(X,t \models\phi)
-
\mathbb{P}((\tilde X,t) \models\phi)
\right|
\leq
\mathbb{P}\left(
|S_\phi^X(t)|\leq C_\phi\eta
\right).
\end{align}
\end{corollary}

\section{Topology and Convergence of TLOPs}
\label{sec_9}
The previous sections established measurability, pathwise regularity, and perturbation properties of TLOPs. We now study the topological structure of the space of TLOPs and investigate the continuity of the mapping
that maps a temporal formula to the TLOP induced by the stochastic process $X$. Particularly, we are interested in establishing the conditions that render the mapping $\mathcal{R}_{X}: \phi \rightarrow S^{\phi}_X$ non-expansive meaning that perturbations in formula $\phi$ do not lead to big deviations in the associated temporal robustness $S^{\phi}_X$. Towards this, we are interested in the following main question:
\begin{problem}
If a sequence of temporal formulas $\{\phi_n\}_{n\in\mathbb{N}}$
converges to a limiting formula $\phi$, under what conditions does the corresponding sequence of TLOPs
$\{S_{\phi_n}^{X}\}_{n\in\mathbb{N}}$ converge to a limiting TLOP $S_{\phi}^{X}$?
\end{problem}

Throughout this section, we restrict attention to the continuous-time case
$\mathbb{T}=\mathbb{R}_{\geq 0}$. Unless stated otherwise, we assume that
$X=\{X_t\}_{t\geq 0}$ has continuous sample paths almost surely and that
the assumptions guaranteeing continuity of the TLOPs established previously
are satisfied.
\subsection{Topology on the Observable-Process Space}
We first introduce a natural topology on the space in which the TLOPs take their values. Fix a finite observation horizon $T>0$. Let $\mathcal{C}_{T} := C([0,T],\mathbb{R})$ denote the space of continuous real-valued functions on $[0,T]$. For $z\in\mathcal{C}_{T}$, define the supremum norm
\begin{align}
    \label{eqn_sup_function_norm}
    \|z\|_{\infty,T}
    :=
    \sup_{t\in[0,T]}|z(t)|.
\end{align}
The choice of the supremum norm is natural for temporal-logical
observables since the quantitative semantics of conjunction, disjunction, always, and eventually are generated by $\min, \max, \inf$, and $\sup$ operations, respectively.
\begin{definition}
\label{def:observable_process_space}
For a fixed finite $T>0$, the observable-process space, denoted by $\mathcal{O}_{T}$ is defined as the space of all random continuous real-valued trajectories on the time interval \([0,T]\), where processes that are equal almost surely are regarded as the same. That is,
\begin{align}
\label{eqn_observable_process_space}
    \mathcal{O}_{T}
    :=
    L^{0}\!\left(
    \Omega,\mathcal{F},\mathbb{P};
    C([0,T],\mathbb{R})
    \right).
\end{align}
\end{definition}
\noindent Under the pathwise regularity assumptions, the restriction of every TLOP $S_\phi^X$ to $[0,T]$ belongs to
$\mathcal{O}_T$. 
\begin{definition}
For $S_{\phi_{1}}^X, S_{\phi_{2}}^X \in\mathcal{O}_{T}$, we define a metric that induces convergence in probability with respect to the
uniform norm as
\begin{align}
    d_{\mathcal{O}_{T}}(S_{\phi_{1}}^X, S_{\phi_{2}}^X)
    :=
    \mathbb{E}
    \left[
    \min \left\{1, 
    \left\| S_{\phi_{1}}^X - S_{\phi_{2}}^X \right\|_{\infty,T}
    \right\}
    \right].
\end{align}
\end{definition}
\begin{proposition}
\label{prop:observable_metric_convergence}
Let $\{S_{\phi_{n}}^X\}_{n\in\mathbb{N}}\subset\mathcal{O}_{T}$ and
$S_{\phi}^X \in \mathcal{O}_{T}$. Then, convergence in $\mathcal{O}_{T}$ corresponds to uniform convergence on $[0,T]$ in probability. That is, 
\begin{align}
    d_{\mathcal{O}_{T}}\left(S_{\phi_{n}}^X, S_{\phi}^X \right)
    \rightarrow 0
    \iff
    \left\|S_{\phi_{n}}^X - S_{\phi}^X \right\|_{\infty,T}
    \xrightarrow{\mathbb{P}}0.
\end{align}
\end{proposition}
\begin{proof}
Let $Z_n := \min \left\{1, \,  \left\|S_{\phi_{n}}^X - S_{\phi}^X \right\|_{\infty,T}\right\}$. Suppose first that $\mathbb{E}[Z_n]\rightarrow 0$. Then, $\forall \epsilon\in(0,1)$, Markov's inequality gives
\begin{align*}
    \mathbb{P}
    \left(
    \left\|S_{\phi_{n}}^X - S_{\phi}^X \right\|_{\infty,T}
    >
    \epsilon
    \right)
    =
    \mathbb{P}(Z_n>\epsilon)
    \leq
    \frac{\mathbb{E}[Z_n]}{\epsilon}
    \rightarrow
    0.
\end{align*}
This means that the probability that the two entire TLOP curves differ by more than \(\epsilon\) anywhere on \([0,T]\) tends to zero. Conversely, let $\left\|S_{\phi_{n}}^X - S_{\phi}^X \right\|_{\infty,T} \xrightarrow{\mathbb{P}}0$. Then $Z_n\xrightarrow{\mathbb{P}}0$. Since $0\leq Z_n\leq 1 \implies \mathbb{E}[Z_n]\rightarrow0$. Hence $d_{\mathcal{O}_{T}}\left(S_{\phi_{n}}^X, S_{\phi}^X \right)\rightarrow0$.
\end{proof}

\subsection{A Semantic Pseudo-metric on Temporal Formulas}
We now introduce a semantic notion of distance between temporal formulas. Metrics between temporal logic specifications have previously been
investigated in the context of Signal Temporal Logic (STL). In
particular, authors in \citep{Madsen2018Metrics} define metrics between STL formulae through distances between their bounded-time satisfaction languages, including a Pompeiu--Hausdorff metric and a symmetric-difference metric. In contrast, our objective here is to introduce a robustness-semantic pseudometric that is directly compatible with the mapping $\phi\mapsto S_\phi^X$. Fix a common temporal-horizon bound $H^{\star}<\infty$ and also fix $T>0$. Let $\Phi_{H^{\star}}$ denote the class of formulas satisfying $H_{\phi}\leq H^{\star}$. Let $\mathcal{K} \subset C([0,T+H^{\star}], \mathcal{X})$ be a collection of admissible trajectories. 
\begin{definition} \label{def:robustness_semantic_metric}
For any $\varphi, \psi \in \Phi_{H^\star}$, define the robustness-semantic distance and its bounded counterpart as
\begin{align}
\label{eqn_distance_btw_formulae}
d_{\rho,T,\mathcal K}(\varphi,\psi)
&:=
\sup_{x\in\mathcal K}
\sup_{t\in[0,T]}
|\rho_\varphi(x,t)-\rho_\psi(x,t)|, \\
\overline d_{\rho,T,\mathcal K}(\varphi,\psi)
&:=
\min\{1,d_{\rho,T,\mathcal K}(\varphi,\psi)\}. \label{eqn_finitepseudodistance_btw_formulae}
\end{align}
\end{definition}
The mapping $d_{\rho,T,\mathcal K}$ defines an extended pseudometric on $\Phi_{H^\star}$, taking values in $[0,\infty]$. The bounded pseudometric  $\overline d_{\rho,T,\mathcal K}(\varphi,\psi)$ takes values
in $[0,1]$ and induces the same notion of convergence as $d_{\rho,T,\mathcal K}$. The distance between formulas may vanish even when the formulas are syntactically different. This motivates the following
semantic equivalence relation.
\begin{definition}
For $\varphi,\psi\in\Phi_{H^\star}$, we define the semantic equivalence relation as
\begin{equation*}
\varphi\sim_\rho\psi
\iff
\rho_\varphi(x,t)=\rho_\psi(x,t),
\quad
\forall x\in\mathcal K,
\forall t\in[0,T].
\end{equation*}
\end{definition}
\noindent Similar quotient constructions based on semantic equivalence have been considered for metrics on STL formulae
in~\citep{Madsen2018Metrics}. Both the bounded pseudometric
$\overline d_{\rho,T,\mathcal K}$ and the original distance
$d_{\rho,T,\mathcal K}$ induce a well-defined metric and an extended metric respectively on the same quotient space $\Phi_{H^\star}/ \sim_\rho$. The bounded semantic pseudometric preserves the convergence structure induced by the original extended pseudometric. Indeed,
\begin{equation*}
\overline d_{\rho,T,\mathcal K}
(\varphi_n,\varphi)\to0
\iff
d_{\rho,T,\mathcal K}
(\varphi_n,\varphi)\to0.
\end{equation*}
The following theorem shows that whenever $X|_{[0,T+H^\star]}\in\mathcal K$ almost surely, the metric $d_{\mathcal O_T}
(S_\varphi^X,S_\psi^X)$ on the observable-process space never exceeds the $\overline d_{\rho,T,\mathcal K}
(\varphi,\psi)$ meaning that the convergence of temporal formulas
in the bounded semantic pseudometric
implies convergence of their induced
TLOPs in $\mathcal O_T$.
\begin{theorem}
\label{thm:nonexpansive_formula_tlop}
Given $\phi,\psi\in\Phi_{H^{\star}}$, let $X(\omega)|_{[0,T+H^{\star}]}\in\mathcal{K}$ almost surely. Then, almost surely, the following holds:
\begin{align}
    \left\|
    S_{\phi}^{X}
    -
    S_{\psi}^{X}
    \right\|_{\infty,T}
    &\leq
    d_{\rho,T,\mathcal{K}}(\phi,\psi), \\
    d_{\mathcal{O}_{T}}
    \left(
    S_{\phi}^{X},
    S_{\psi}^{X}
    \right)
    &\leq
    \overline d_{\rho,T,\mathcal K} (\varphi,\psi).
\end{align}
Further, if $d_{\rho,T,\mathcal{K}}(\phi_n, \phi) \rightarrow0$, then $\left\| S_{\phi_n}^{X} -
S_{\phi}^{X} \right\|_{\infty,T}
\rightarrow0$ almost surely, and therefore $S_{\phi_n}^{X} \rightarrow S_{\phi}^{X} \in \mathcal{O}_{T}$.
\end{theorem}

\begin{proof}
For almost every $\omega\in\Omega$,
\begin{align*}
    \left\|
    S_{\phi}^{X}
    -
    S_{\psi}^{X}
    \right\|_{\infty,T}
    =
    \sup_{t\in[0,T]}
    \left|
    \rho_{\phi}(X(\omega),t)
    -
    \rho_{\psi}(X(\omega),t)
    \right|
    \leq
    \sup_{x\in\mathcal{K}}
    \sup_{t\in[0,T]}
    \left|
    \rho_{\phi}(x,t)
    -
    \rho_{\psi}(x,t)
    \right|
    =
    d_{\rho,T,\mathcal{K}}(\phi,\psi).
\end{align*}
Taking the expectation of the truncated supremum norm gives the second
claim.
\end{proof}
\subsection{Structural-Parametric Convergence of Temporal Formulas}
The previous result defines convergence directly through the robustness semantics. Now consider the formula sequence defined by $\phi_n := \Box_{[0, 5]} (X \geq 10+1/n)$. Since as $n \rightarrow \infty, 10+1/n \rightarrow 10$, we infer that $\phi_n \rightarrow \phi:= \Box_{[0, 5]} (X \geq 10)$. Here, the formula structure stays the same while the numerical parameter converges. We now introduce a more constructive notion of convergence based on the parameters appearing in temporal formulas.

\begin{definition}
\label{def:structural_parametric_convergence}
Let $\{\phi_n\}_{n\in\mathbb{N}}$ and $\phi$ be temporal formulas.
We write $\phi_n \xrightarrow{\mathrm{sp}} \phi$ if the following conditions hold.
\begin{enumerate}
    \item Both $\phi_n$ and $\phi$ have the same finite
    syntactic tree, consisting of  operators from \eqref{eqn_formula_semantics}.
    \item Every sequence of atomic predicates $\mu_n$ and limiting predicate $\mu$, with predicate functions $h_{\mu_n}$ and $h_{\mu}$, respectively satisfy $h_{\mu_n} \rightarrow h_{\mu}$ locally uniformly on $\mathcal{X}$.
    \item Every sequence of temporal intervals $I_n=[a_n,b_n]$ and limiting interval $I=[a,b]$ satisfy $a_n \rightarrow a, b_n \rightarrow b$.
    \item The temporal horizons are uniformly bounded meaning that $\sup_{n\in\mathbb{N}} H_{\phi_n} <\infty$.
\end{enumerate}
\end{definition}
\noindent For compact intervals, convergence of the interval endpoints is equivalent to convergence with respect to the Hausdorff metric meaning that $d_{\mathrm H}([a_n,b_n],[a,b]) = \max \left\{ |a_n-a|, |b_n-b| \right\}$. We now establish the main continuity result.
\begin{theorem}
\label{thm:parametric_formula_convergence}
Assume that every atomic predicate function appearing in $\phi_n$ and $\phi$ is continuous.
Let $\phi_n \xrightarrow{\mathrm{sp}} \phi$. Then, for every finite $T>0$, we see that
\begin{align}
    \sup_{t\in[0,T]}
    \left|
    S_{\phi_n}^{X}(t)
    -
    S_{\phi}^{X}(t)
    \right|
    \rightarrow0
\end{align}
almost surely. Consequently, $S_{\phi_n}^{X} \rightarrow
S_{\phi}^{X} \in \mathcal{O}_{T}$ for every finite $T>0$.
\end{theorem}

\begin{proof}
Let $\Omega_c\subseteq\Omega$ denote the event on which the sample path
$t\mapsto X_t(\omega)$ is continuous. By assumption, $\mathbb{P}(\Omega_c)=1$. Fix $\omega\in\Omega_c$. Since every temporal formula is obtained by finitely many applications of the grammar rules, we will prove the result by structural induction. \\
Case 1: Atomic predicates. Suppose $\phi_n=\mu_n$ and $\phi=\mu$. Then, $S_{\mu_n}^{X}(t)(\omega)
= h_{\mu_n}(X_t(\omega)), S_{\mu}^{X}(t)(\omega)
= h_{\mu}(X_t(\omega))$. Since $t\mapsto X_t(\omega)$ is continuous, the set $K_{\omega,T} := \left\{ X_t(\omega) \mid t\in[0,T] \right\}$ is compact. Local uniform convergence of $h_{\mu_n}$ to $h_{\mu}$ gives
\begin{align*}
    \sup_{t\in[0,T]}
    \left|
    h_{\mu_n}(X_t(\omega))
    -
    h_{\mu}(X_t(\omega))
    \right|
    \leq
    \sup_{x\in K_{\omega,T}}
    |h_{\mu_n}(x)-h_{\mu}(x)|
    \rightarrow0.
\end{align*}
Case 2: Negation. Suppose $\phi_n = \neg\psi_n, \phi=\neg\psi$. Then
$\left|
    S_{\phi_n}^{X}(t)
    -
    S_{\phi}^{X}(t)
    \right|
    =
    \left|
    S_{\psi_n}^{X}(t)
    -
    S_{\psi}^{X}(t)
\right|$. The result follows from the induction hypothesis. \\
Case 3: Conjunction and disjunction. Suppose $\phi_n = \psi_{1,n}\wedge\psi_{2,n}, \phi = \psi_1\wedge\psi_2$. Using the elementary inequality
\begin{align*}
    |\min(a,b)-\min(c,d)|
    \leq
    \max
    \left\{
    |a-c|,
    |b-d|
    \right\},
\end{align*}
we obtain
\begin{align*}
    \sup_{t\in[0,T]}
    \left|
    S_{\phi_n}^{X}(t)
    -
    S_{\phi}^{X}(t)
    \right|
    \leq
    \max
    \left\{
    \sup_{t\in[0,T]}
    \left|
    S_{\psi_{1,n}}^{X}(t)
    -
    S_{\psi_1}^{X}(t)
    \right|,
    \sup_{t\in[0,T]}
    \left|
    S_{\psi_{2,n}}^{X}(t)
    -
    S_{\psi_2}^{X}(t)
    \right|
    \right\}.
\end{align*}
The right-hand side converges to zero by the induction hypothesis.
The case of disjunction follows analogously using the corresponding
inequality for the maximum. \\
Case 4: Always operator. Suppose $\phi_n = \Box_{I_n}\psi_n, \phi = \Box_I\psi$, where $I_n=[a_n,b_n], \,I=[a,b]$. Since $b_n\rightarrow b$, there exists a finite constant $B>0$ such that $b_n\leq B, b\leq B$ for all sufficiently large $n$. Define $r_n(s) := S_{\psi_n}^{X}(s)(\omega), r(s) := S_{\psi}^{X}(s)(\omega)$. By the induction hypothesis, $\|r_n-r\|_{\infty,T+B} \rightarrow0$. Moreover, $r$ is continuous on the compact interval $[0,T+B]$ and is therefore uniformly continuous. Define the modulus of continuity of $r$ on this interval by
\begin{align*}
    \omega_r(\delta)
    :=
    \sup
    \left\{
    |r(u)-r(v)|
    \,\middle|\,
    u,v\in[0,T+B],
    \ |u-v|\leq\delta
    \right\}.
\end{align*}
Then, $\omega_r(\delta)\rightarrow0$
as $\delta\rightarrow0$. For every $t\in[0,T]$, we see that $S_{\phi_n}^{X}(t)(\omega) = \inf_{s \in t+I_n}r_n(s)$, and $S_{\phi}^{X}(t)(\omega) = \inf_{s\in t+I}r(s)$. Hence,
\begin{align*}
    \left|
    S_{\phi_n}^{X}(t)(\omega)
    -
    S_{\phi}^{X}(t)(\omega)
    \right|
    &\leq
    \left|
    \inf_{s\in t+I_n}r_n(s)
    -
    \inf_{s\in t+I_n}r(s)
    \right|
    +
    \left|
    \inf_{s\in t+I_n}r(s)
    -
    \inf_{s\in t+I}r(s)
    \right| \\
    &\leq
    \|r_n-r\|_{\infty,T+B}
    +
    \omega_r
    \left(
    d_{\mathrm H}(I_n,I)
    \right).
\end{align*}
Both terms in the right hand side converge to zero, and hence the result follows. \\
Case 5: Eventually operator. Suppose $\phi_n =
\Diamond_{I_n}\psi_n, \phi =
\Diamond_I\psi$. The proof is identical to the previous case, using $S_{\phi_n}^{X}(t) = \sup_{s\in t+I_n} S_{\psi_n}^{X}(s)$ and the inequality $\left| \sup f-\sup g \right| \leq \sup|f-g|$. 
\end{proof}

\subsection{Weak Convergence of TLOP Distributions}
Since a TLOP is a random element of $C([0,T],\mathbb{R})$, it induces a
probability measure on this function space. For a TLOP $S_{\phi}^{X}$, denote its law by $\mathcal{L}(S_{\phi}^{X}) \in \mathcal{P}(C([0,T], \mathbb{R}))$. The following corollary says that the whole statistical behaviour of the temporal-logical process converges.
\begin{corollary}
\label{cor:weak_convergence_tlop_laws}
Suppose $S_{\phi_n}^{X}\rightarrow S_{\phi}^{X}$ in $\mathcal{O}_{T}$. Then $\mathcal{L}(S_{\phi_n}^{X})
\Rightarrow \mathcal{L}(S_{\phi}^{X})$, where $\Rightarrow$ denotes weak convergence of probability measures on $C([0,T],\mathbb{R})$. Consequently, for every bounded continuous functional $F \in C([0,T],\mathbb{R}) \rightarrow \mathbb{R}$, we have
\begin{align}
\label{eqn_convergence_expected_fn_of_TLOP}
    \mathbb{E}
    \left[
    F(S_{\phi_n}^{X})
    \right]
    \rightarrow
    \mathbb{E}
    \left[
    F(S_{\phi}^{X})
    \right].
\end{align}
\end{corollary}

\begin{proof}
By the definition of the metric $d_{\mathcal{O}_T}$ and
Proposition~\ref{prop:observable_metric_convergence}, the convergence $S_{\phi_n}^{X}\to S_{\phi}^{X}$ in $\mathcal{O}_T$ is equivalent to $\left\|
S_{\phi_n}^{X} - S_{\phi}^{X} \right\|_{\infty,T} \xrightarrow{\mathbb{P}} 0$. Hence, when regarded as random elements of the metric space $C([0,T],\mathbb{R})$ equipped with the supremum norm,
$S_{\phi_n}^{X}$ converges in probability to $S_{\phi}^{X}$. Convergence in probability implies convergence in distribution, and therefore $\mathcal{L}(S_{\phi_n}^{X}) \Rightarrow \mathcal{L}(S_{\phi}^{X})$. Since $F:C([0,T],\mathbb{R})\to\mathbb{R}$ is bounded and continuous, the convergence of expectations follows from the definition of weak convergence.
\end{proof}

\subsection{Closedness of Satisfaction and Loss of Strict Robustness}
We now investigate whether temporal satisfaction is preserved under convergence of TLOPs. In particular, we distinguish between the non-negativity of quantitative robustness, strictly positive robustness, and Boolean satisfaction. Although non-negative robustness is preserved under uniform
convergence, this does not necessarily imply preservation of Boolean satisfaction. As discussed before, zero robustness represents a logical boundary at which Boolean satisfaction cannot always be determined from
the robustness sign alone. We first establish the closedness of non-negative TLOP trajectories under uniform convergence.
\begin{proposition}
Let $S_{\varphi_n}^X\to S_\varphi^X$ uniformly on $[0,T]$ for some finite $T>0$ almost surely. Assume further that $S_{\varphi_n}^X(t)\geq 0, \forall n \in \mathbb{N}, \forall t\in[0,T]$, almost surely, with a common probability-one event. Then the limiting TLOP satisfies $S_\varphi^X(t)\geq 0, \forall t\in[0,T]$, almost surely.
\end{proposition}
\begin{proof}
Let $\Omega_0$ be the intersection of the probability-one event on which uniform convergence holds and the
probability-one event on which all approximating TLOPs are non-negative $\forall t\in[0,T]$. Fix $\omega\in\Omega_0$. Uniform convergence implies 
\begin{equation*}
\sup_{t\in[0,T]}
\left|
S_{\varphi_n}^X(t)(\omega)
-
S_\varphi^X(t)(\omega)
\right|
\longrightarrow 0.
\end{equation*}
Consequently, $\forall t\in[0,T]$, we see that $S_{\varphi_n}^X(t)(\omega) \longrightarrow S_\varphi^X(t)(\omega)$. Since $S_{\varphi_n}^X(t)(\omega)\geq 0, \forall n$, the closedness of $[0,\infty)$ implies $S_\varphi^X(t)(\omega)\geq 0$. Since $t$ and $\omega$ were arbitrary, the result follows.
\end{proof}
The preceding proposition establishes closedness of non-negative robustness. However, strictly positive robustness need not be preserved. More importantly,
Boolean satisfaction itself may fail to be preserved when the limiting robustness reaches zero. The following proposition demonstrates both phenomena.
\begin{proposition}
There exist a stochastic process $X$, a sequence of temporal formulas $\{\varphi_n\}_{n\in\mathbb{N}}$, and a limiting formula $\varphi$ such that $S_{\varphi_n}^X(t)>0, \forall n\in\mathbb{N}$, and $S_{\varphi_n}^X\longrightarrow S_\varphi^X=0$ uniformly on every compact time interval, while $(X,t)\models\varphi_n, \forall n\in\mathbb{N}$, but $(X,t) \not\models \varphi$. 
\end{proposition}

\begin{proof}
Consider the constant stochastic process $X_t(\omega) = e, \forall t\geq 0, \forall\omega\in\Omega$. Let $\{q_n\}_{n\in\mathbb{N}}\subset\mathbb{Q}$ satisfy $q_n>e$ for every $n$ and $q_n\to e$. Define the temporal formulas
\begin{align*}
\varphi_n
:=\Box_{[0,T]}\neg(X\geq q_n), \quad 
\varphi
:=\Box_{[0,T]}\neg(X\geq e).
\end{align*}
Clearly, $\varphi_n$ converges structurally and parametrically to $\varphi$ according to Definition \ref{def:structural_parametric_convergence}. Since $X_t=e<q_n$, the atomic predicate $(X\geq q_n)$ is false at every time, and hence $(X,t)\models\varphi_n, \forall n\in\mathbb{N}, \forall t\geq0$. The corresponding TLOP is
\begin{align*}
S_{\varphi_n}^X(t)
=
\inf_{s\in[t,t+T]}(q_n-X_s) 
=q_n-e>0.
\end{align*}
Since $q_n\to e$, we obtain $S_{\varphi_n}^X(t)\longrightarrow 0
=S_\varphi^X(t)$. This convergence is uniform on every compact time interval since the observables are constant
in time. Given the limiting formula $\varphi$ and since $X_t=e$, the predicate $(X\geq e)$ is true at every time, and therefore $(X,t) \not\models \varphi, \forall t\geq0$. Thus, the approximating formulas are all satisfied
with strictly positive robustness, while their limiting formula is violated and has zero robustness.
\end{proof}
The preceding proposition shows that the loss of strict robustness may also lead to a change in Boolean satisfaction. Such changes can occur at the zero-robustness boundary. Strictly positive robustness is, however,
preserved when a common positive margin exists.
\begin{corollary}
Suppose that $S_{\varphi_n}^X\to S_\varphi^X$ uniformly on $[0,T]$ almost surely. Assume further that $\exists \delta>0$ such that $S_{\varphi_n}^X(t)\geq\delta,
\forall n\in\mathbb{N}, \forall t \in [0,T]$ almost surely, with a common probability-one event. Then $S_\varphi^X(t)\geq\delta>0, \forall t\in[0,T]$ almost surely. Consequently, $(X,t)\models\varphi, \forall t\in[0,T]$
almost surely.
\end{corollary}

\begin{proof}
On the common probability-one event of uniform convergence and the uniform positive-margin condition, taking the limit in $S_{\varphi_n}^X(t)\geq\delta$ yields $S_\varphi^X(t)\geq\delta>0, \forall t\in[0,T]$. The Boolean satisfaction conclusion follows from
the sign-soundness property of quantitative robustness established in Section 2.2.
\end{proof}

\subsection{Convergence of Satisfaction Probabilities}
Convergence of TLOPs does not automatically imply convergence of the associated probability of satisfaction because the mapping $z\longmapsto\mathbf{1}_{{z\geq0\}}}$ is discontinuous at $z=0$. So, the probability mass at the logical
boundary must be treated separately.
\begin{theorem}
\label{thm:satisfaction_probability_convergence}
Fix $t\geq 0$ and let $S_{\phi_n}^{X}(t) \xrightarrow{\mathbb{P}} S_{\phi}^{X}(t)$ and that $\mathbb{P}(S_{\phi}^{X}(t)=0)=0$. Then, $\mathbb{P}\big(X,t\models\phi_n\big)
\longrightarrow \mathbb{P}\big(X,t\models\phi\big)$.
\end{theorem}
\begin{proof}
Since $S_{\phi_n}^{X}(t)\xrightarrow{\mathbb{P}} S_{\phi}^{X}(t)$, we also have $S_{\phi_n}^{X}(t)\Rightarrow S_{\phi}^{X}(t)$. Because $\mathbb{P}(S_{\phi}^{X}(t)=0)=0$, the point $0$ is a continuity point of the distribution function of $S_{\phi}^{X}(t)$. Hence, $\mathbb{P}(S_{\phi_n}^{X}(t)>0) \longrightarrow \mathbb{P}(S_{\phi}^{X}(t)>0)$ and $\mathbb{P}(S_{\phi_n}^{X}(t)\geq 0) \longrightarrow \mathbb{P}(S_{\phi}^{X}(t)\geq 0)$. From the sign-soundness property of the robustness semantics,
\[
\{S_{\phi_n}^{X}(t)>0\}
\subseteq
\{
\omega \mid X(\omega),t\models\phi_n
\}
\subseteq
\{S_{\phi_n}^{X}(t)\geq 0\}.
\]
Taking probabilities, we get
\[
\mathbb{P}(S_{\phi_n}^{X}(t)>0)
\leq
\mathbb{P}\big(X,t\models\phi_n\big)
\leq
\mathbb{P}(S_{\phi_n}^{X}(t)\geq0).
\]
Similarly, since $\mathbb{P}(S_{\phi}^{X}(t)=0)=0$,
\[
\mathbb{P}\big(X,t \models\phi\big)
=
\mathbb{P}(S_{\phi}^{X}(t)>0)
=
\mathbb{P}(S_{\phi}^{X}(t)\geq0).
\]
Both the lower and upper bounds for
$\mathbb{P}(X,t\models\phi_n)$ therefore converge to $\mathbb{P}(X,t\models\phi)$. The result follows from the squeeze theorem.
\end{proof}
The assumption $\mathbb{P} \left( S_{\phi}^{X}(t)=0 \right) = 0$ is crucial as the logical boundary $S_{\phi}^{X}(t)=0$
is precisely the point at which arbitrarily small perturbations of the TLOP can change the sign of the robustness. This observation is closely related to the probability-level robustness
result established earlier, where the quantity $\mathbb{P} \left( |S_{\phi}^{X}(t)| \leq\epsilon \right)$ measures the probability mass concentrated near the logical boundary.
\begin{remark}
The preceding convergence results naturally motivate the study of the completion of the temporal-formula space under a suitable semantic metric. An interesting open question is whether such a completion contains limit
objects whose induced TLOPs cannot be represented by any finite temporal formula in $\Phi$. Such objects, if they exist, may be interpreted as generalized temporal specifications.
\end{remark}

\section{Large Deviations of Temporal-Logical Observable Processes}
\label{sec_10}

The preceding sections establish that a temporal-logical observable
process inherits useful measurability, regularity, perturbation, and
convergence properties from the underlying stochastic process. We now
consider a different question. Suppose that a family of stochastic
processes depends on a small noise parameter $\epsilon>0$, and that
a temporal event becomes increasingly unlikely as $\epsilon\to0$.
Then, we are interested in the probability of this rare temporal event being characterized directly through the corresponding TLOP. The natural framework for answering this question is the theory of
large deviations. We first illustrate the main idea using a simple
example.

\subsection{A Motivating Example: Rare Target Attainment}
Given $\epsilon>0$, and $x_0<K$, let the small-noise Brownian process be given by
\begin{align}
\label{eq:small_noise_BM}
X_t^{\epsilon}
= x_0+\sqrt{\epsilon}\,W_t, \quad t\in[0,T],  
\end{align}
where $W=\{W_t\}_{t\geq0}$ is a standard Brownian motion. Consider the temporal specification
\begin{align}
\label{eq:rare_event_formula}
\phi
=
\lozenge_{[0,T]}(X\geq K),
\end{align}
requiring the process to reach the level $K$ within the time interval $[0,T]$. From Theorem \ref{thm_reflection_principle}, we know that 
\begin{align}
\label{eq:exact_small_noise_probability}
\mathbb{P}\big((X^\epsilon,0)\models\phi\big)
&=
\mathbb{P}\left(
\sup_{0\leq s\leq T}
X_s^\epsilon\geq K
\right)
=
2\left[
1-\Phi\left(
\frac{K-x_0}{\sqrt{\epsilon T}}
\right)
\right], \\
\implies 
\lim_{\epsilon\downarrow0}
\epsilon
\log
\mathbb{P}\big((X^\epsilon,0)\models\phi\big)
&=
-\frac{(K-x_0)^2}{2T},
\label{eq:Brownian_LDP_exponent}
\end{align}
where we used the standard Gaussian-tail asymptotics to get the above implication. More precisely, using the Gaussian-tail asymptotic
\begin{align}
\label{eqn_rare_target_attainment_Gaussian_tail_asymptotics}
1-\Phi_N(z)
&\sim
\frac{1}{z\sqrt{2\pi}}
\exp\left(-\frac{z^2}{2}\right),
\quad \text{as } z\to\infty, \text{ and so}, \\
\mathbb{P}\left(
(X^\varepsilon,0)\models\varphi
\right)
&\sim
\frac{\sqrt{2\varepsilon T/\pi}}{K-x_0}
\exp\left(
-\frac{(K-x_0)^2}{2T\varepsilon}
\right),
\qquad \varepsilon\downarrow 0.
\label{eqn_rare_target_attainment_Gaussian_tail_equivalence_relation}
\end{align}
Thus, the probability of target attainment decays exponentially as $\varepsilon\downarrow 0$, with the logarithmic decay rate given by \eqref{eq:Brownian_LDP_exponent}. The overloaded symbol $\sim$ in \eqref{eqn_rare_target_attainment_Gaussian_tail_asymptotics} and \eqref{eqn_rare_target_attainment_Gaussian_tail_equivalence_relation} denote asymptotic equivalence, meaning that the ratio of the two expressions converges to one. For a continuous trajectory $x:[0,T] \to \mathbb{R}$ starting at $x_0$, define
\begin{align}
I_T(x)
=
\begin{cases}
\displaystyle
\frac{1}{2}\int_0^T |\dot{x}(s)|^2\,ds,
&
x\text{ absolutely continuous and }x(0)=x_0,
\\[2mm]
+\infty,
&
\text{otherwise}.
\end{cases}
\label{eq:Brownian_action}
\end{align}
This is the sample-path rate function associated with
\eqref{eq:small_noise_BM}, as follows from Schilder's theorem; see, for example, \citep{DemboZeitouni2010}. The temporal event in \eqref{eq:rare_event_formula} corresponds to the satisfaction
set
\begin{align}
\label{eqn_A_phi}
\mathcal{A}_\phi
:=
\left\{
x\in C([0,T],\mathbb{R})
\mid
\sup_{0\leq s\leq T}x(s)\geq K
\right\}.
\end{align}
The following proposition established the minimum action required to satisfy this temporal specification, whose satisfaction set is given by \eqref{eqn_A_phi}.
\begin{proposition}
\label{proposition_minimum_action_temporal}
The minimum action required to satisfy the temporal specification $\phi$ defined in \eqref{eq:rare_event_formula} is given by 
\begin{align}
\label{eq:minimum_Brownian_action}
I_\phi^\star
:=
\inf_{x\in\mathcal{A}_\phi}I_T(x)
=
\frac{(K-x_0)^2}{2T}.
\end{align}
\end{proposition}
\begin{proof}
If an admissible path first reaches $K$ at time $\tau\leq T$, then from fundamental theorem of calculus, we know that
\begin{align}
K-x_0
=
\int_0^\tau \dot{x}(s)\,ds.
\end{align}
By the Cauchy--Schwarz inequality,
\begin{align*}
(K-x_0)^2
\leq
\tau
\int_0^\tau |\dot{x}(s)|^2\,ds 
\iff
\frac{(K-x_0)^2}{2\tau}
\leq
I_{T}(x) 
\implies
I_T(x)
\geq
\frac{(K-x_0)^2}{2\tau}
\geq
\frac{(K-x_0)^2}{2T}.
\end{align*}
Any path of $X^{\epsilon}$ reaching $K$ by time $\tau \leq T$ must cost at least $\frac{(K-x_0)^2}{2T}$. Equality in Cauchy Schwarz occus when $f(s) = c g(s)$ for some constant $c \in \mathbb{R}$ with $f(s) = \dot{x}(s)$ and $g(s) = 1$ in this case. This means that equality is attained when $\dot{x}(s) = c$ meaning that velocity must be constant. That is, $x(s) = x_0 + cs$. Since, $x(\tau) = K$, we see that 
\begin{align*}
K 
&= 
x_0 + c\tau 
\iff c = \frac{K - x_0}{\tau}
\end{align*}
The cheapest path reaching $K$ exactly at time $\tau$ is attained by the linear trajectory from $(0, x_0)$ to $(\tau, K)$. Hence,
\begin{align}
x^\star(s)
=
x_0+\frac{K-x_0}{T}s,
\qquad
0\leq s\leq T.
\label{eq:most_likely_Brownian_path}
\end{align}
Consequently, the equality of cauchy Schwarz follows and hence the result \eqref{eq:minimum_Brownian_action} follows too. 
\end{proof}

Comparison of \eqref{eq:Brownian_LDP_exponent} and
\eqref{eq:minimum_Brownian_action} reveals the fundamental principle
that motivates this section that the probability of a rare temporal event decreases exponentially in terms of the 
minimum pathwise cost required to realize that temporal event.
We now formalize this principle for general TLOPs.

\subsection{Large-deviation principle}
We briefly recall the minimum large-deviation machinery required below. For a systematic treatment, see \citep{DemboZeitouni2010}.

\begin{definition}
Let $\mathcal{E}$ be a metric space. A function $I : \mathcal{E} \to [0,\infty]$ is called a rate function if it is lower semicontinuous. In addition, it is called a good rate function if every level set $\{x \in \mathcal{E} \mid I(x) \leq \alpha \}$ is compact $\forall \alpha < \infty$.
\end{definition}

\begin{definition}
Let $\{Y^\epsilon\}_{\epsilon>0} \in \mathcal{E}$ be a family of random elements. The family is said to satisfy a large-deviation principle  (LDP) with speed $\frac{1}{\epsilon}$ and rate function $I$ if, for every Borel set $A \subseteq \mathcal{E}$,
\begin{align}
-\inf_{y\in A^\circ}I(y)
\leq
\liminf_{\epsilon\downarrow0}
\epsilon
\log\mathbb{P}(Y^\epsilon\in A)
\leq
\limsup_{\epsilon\downarrow0}
\epsilon
\log\mathbb{P}(Y^\epsilon\in A)
\leq
-\inf_{y\in\overline{A}}I(y),
\label{eq:LDP_definition}
\end{align}
where $A^\circ$ and $\overline A$ denote the interior and closure of
$A$, respectively. Whenever the infima over $A^\circ$ and $\overline{A}$ agree,
\begin{align*}
\lim_{\epsilon\downarrow0}
\epsilon
\log\mathbb{P}(Y^\epsilon\in A)
=
-\inf_{y\in A}I(y) 
\iff
\mathbb{P}(Y^\epsilon\in A)
\approx
e^{-\frac{1}{\epsilon} \inf_{y\in A}I(y)}
\end{align*}
\end{definition}
\noindent The dominant contribution to the rare-event probability therefore comes from the least costly element of the event. The rate function $I(y)$ captures the exponential cost of
observing the rare outcome $y$. Outcomes with $I(y)=0$ correspond
to typical limiting behavior, whereas larger values of $I(y)$ describe
increasingly unlikely deviations. To understand this phenomenon, consider an ordinary Gaussian random variable $X^{\epsilon} := \sqrt{\epsilon} Z$, where $Z \sim \mathcal{N}(0, 1)$ and $\epsilon > 0$. When $\epsilon$ is small, $X^{\epsilon}$ is concentrated very close to zero because $\operatorname{Var}(X^{\epsilon}) = \epsilon$. Then, given a fixed $a > 0$, we see that
\begin{align*}
\mathbb{P}(X^{\epsilon} \geq a)
= 
\mathbb{P}(\sqrt{\epsilon} Z \geq a)
= 
\mathbb{P}\left(Z \geq \frac{a}{\sqrt{\epsilon}}\right)
\approx
e^{-\frac{a^2}{2 \epsilon}},
\end{align*}
where we used the fact that $\frac{a}{\sqrt{\epsilon}} \rightarrow \infty$ (becomes an extremely unlikely Gaussian-tail event) as $\epsilon \rightarrow 0$. Equivalently, we can write the above probability as
\begin{align*}
\lim_{\epsilon \rightarrow 0} \epsilon 
\log \mathbb{P}(X^{\epsilon} \geq a)
= 
-\frac{a^2}{2}.
\end{align*}
The quantity $I(a) = \frac{a^2}{2}$ is the cost of producing the rare deviation $a$.  
\subsection{The temporal-robustness path map}

Fix $T>0$ and a temporal formula $\phi$ with finite horizon $H_\phi$. Consider the path spaces $\mathcal{C}_{T+H_\phi}^{d}
:=
C([0,T+H_\phi],\mathbb{R}^d), \mathcal{C}_{T}
:=
C([0,T],\mathbb{R})$, both equipped with the uniform norm. Define the temporal-robustness map $\mathcal{R}_{\phi,T}
: \mathcal{C}_{T+H_\phi}^{d} \longrightarrow \mathcal{C}_{T}$ by
\begin{align}
\big(\mathcal{R}_{\phi,T}(x)\big)(t)
:=
\rho_\phi(x,t),
\qquad
t\in[0,T].
\label{eq:robustness_path_map_definition}
\end{align}
Thus, if $X$ is an underlying stochastic process, then $\mathcal{R}_{\phi,T}(X)
= \left\{ S_\phi^X(t)
\right\}_{t\in[0,T]}$. The pathwise perturbation result established earlier immediately yields
continuity of this map.

\begin{proposition}\label{prop:robustness_map_continuity}
Suppose that every temporal interval appearing in $\phi$ is compact and that every atomic robustness function is globally Lipschitz. Then, $\mathcal{R}_{\phi,T}$ is continuous meaning that $\exists C_\phi>0$ such that, $\forall x, y\in\mathcal{C}_{T+H_\phi}^{d}$,
\begin{align}
\left\|
\mathcal{R}_{\phi,T}(x)
-
\mathcal{R}_{\phi,T}(y)
\right\|_{\infty,T}
\leq
C_\phi
\left\|
x-y
\right\|_{\infty,T+H_\phi}.
\label{eq:robustness_map_Lipschitz}
\end{align}
\end{proposition}

\begin{proof}
The estimate follows from the same structural argument used in the
pathwise perturbation Theorem \ref{thm_pathwise_robustness}. Atomic predicates inherit their
Lipschitz constants, while negation, minimum, maximum, temporal
infimum, and temporal supremum are nonexpansive under the uniform norm.
Taking the supremum over $t\in[0,T]$ gives
\eqref{eq:robustness_map_Lipschitz}.
\end{proof}

\subsection{A large-deviation principle for TLOPs}
The continuity property established in Proposition \ref{prop:robustness_map_continuity} is the key connection between large deviations of the original stochastic process and large deviations of its induced
TLOP. We can now state the main result of this section.

\begin{theorem}
\label{thm:TLOP_LDP}
Let $\{X^\epsilon\}_{\epsilon>0}$ be a family of stochastic processes whose sample paths belong almost surely to
$\mathcal{C}_{T+H_\phi}^{d}$ and $\{X^\epsilon\}_{\epsilon>0}$ satisfies an LDP with speed
$\epsilon^{-1}$ and good rate function $I: \mathcal{C}_{T+H_\phi}^{d}
\longrightarrow[0,\infty]$. Assume that Proposition~\ref{prop:robustness_map_continuity} hold. Then the family of TLOPs, $S_\phi^{X^\epsilon} =
\mathcal{R}_{\phi,T}(X^\epsilon)$ satisfies an LDP on $\mathcal{C}_{T}$ with speed $\epsilon^{-1}$ and good rate function\footnote{infimum over the empty set is defined to be $+\infty$}
\begin{align}
J_{\phi,T}(z)
=
\inf_{\substack{
x\in\mathcal{C}_{T+H_\phi}^{d}\\
\mathcal{R}_{\phi,T}(x)=z
}}
I(x),
\qquad
z\in\mathcal{C}_{T}.
\label{eq:TLOP_rate_function}
\end{align}
\end{theorem}

\begin{proof}
By Proposition~\ref{prop:robustness_map_continuity}, $\mathcal{R}_{\phi,T}:
\mathcal{C}_{T+H_\phi}^{d} \to \mathcal{C}_{T}$ is continuous. Since $S_\phi^{X^\epsilon} =
\mathcal{R}_{\phi,T}(X^\epsilon)$, the result follows directly from the contraction principle for large deviations \citep{DemboZeitouni2010}, which also implies that $J_{\phi,T}$ is a good rate function.
\end{proof}

Theorem~\ref{thm:TLOP_LDP} provides a process-level interpretation of
rare temporal behavior. For a candidate logical trajectory
$z\in\mathcal{C}_T$, the quantity $J_{\phi,T}(z)$ is the minimum
large-deviation cost among all physical trajectories whose temporal
robustness trajectory is exactly $z$.

\begin{remark}
The large-deviation principle assumed for $X^\epsilon$ is not tied
to a particular stochastic model. Theorem~\ref{thm:TLOP_LDP} therefore
applies whenever a suitable path-space LDP is available. The role of
the temporal-logic framework is to transfer the known stochastic
large-deviation structure of $X^\epsilon$ through the continuous
map $\mathcal{R}_{\phi,T}$ to obtain a corresponding large-deviation
description of the entire TLOP.
\end{remark}

\subsection{Large deviations of fixed-time temporal robustness}

The process-level result immediately yields a scalar large-deviation
principle for the temporal robustness observed at any fixed time.

\begin{corollary}
\label{cor:fixed_time_robustness_LDP}
Fix $t\in[0,T]$. Then $\left\{
S_\phi^{X^\epsilon}(t)
\right\}_{\epsilon>0}$ satisfies an LDP on $\mathbb{R}$ with rate function
\begin{align}
J_{\phi,t}(r)
=
\inf_{\substack{
x\in\mathcal{C}_{T+H_\phi}^{d}\\
\rho_\phi(x,t)=r
}}
I(x).
\label{eq:fixed_time_rate_function}
\end{align}
\end{corollary}

\begin{proof}
The evaluation mapping $e_t:\mathcal{C}_T\to\mathbb{R}$ given by $e_t(z)=z(t)$ is continuous under the uniform norm. Applying the contraction
principle to Theorem~\ref{thm:TLOP_LDP} gives the result.
\end{proof}

This result also provides asymptotic estimates for temporal
satisfaction and violation probabilities. Because zero robustness is a
logical boundary, it is useful to distinguish strict violation from
the boundary event. Define
\begin{align}
I_{\phi,t}^{<0}
&:=
\inf_{\substack{x\in\mathcal{C}_{T+H_\phi}^{d}\\
\rho_\phi(x,t)<0}}
I(x),
\label{eq:strict_violation_cost}
\\
I_{\phi,t}^{\leq0}
&:=
\inf_{\substack{x\in\mathcal{C}_{T+H_\phi}^{d}\\
\rho_\phi(x,t)\leq0}}
I(x).
\label{eq:closed_violation_cost}
\end{align}

\begin{corollary}
\label{cor:rare_violation}
The following inequality holds due to Theorem~\ref{thm:TLOP_LDP},
\begin{align}
-I_{\phi,t}^{<0}
\leq
\liminf_{\epsilon\downarrow0}
\epsilon
\log
\mathbb{P}\big((X^\epsilon,t)\not\models\phi\big)
\leq
\limsup_{\epsilon\downarrow0}
\epsilon
\log
\mathbb{P}\big((X^\epsilon,t)\not\models\phi\big)
\leq
-I_{\phi,t}^{\leq0}.
\label{eq:rare_violation_bounds}
\end{align}
Consequently, if
\begin{equation}
\label{eq:boundary_rate_condition}
I_{\varphi,t}^{<0}
=
I_{\varphi,t}^{\leq 0},
\end{equation}
we define the temporal-logic failure rate as $I_{\varphi,t}^{\mathrm{fail}} := I_{\varphi,t}^{<0} = I_{\varphi,t}^{\leq 0}$. Then, the large-deviation lower and upper bounds in
\eqref{eq:rare_violation_bounds} coincide, yielding
\begin{equation}
\label{eq:exact_failure_exponent}
\lim_{\varepsilon\downarrow 0}
\varepsilon
\log
\mathbb{P}
\left(
(X^\varepsilon,t)\not\models\varphi
\right)
=
-I_{\varphi,t}^{\mathrm{fail}}.
\end{equation}
\end{corollary}

\begin{proof}
We know from the quantitative semantics that,
\begin{align}
\left\{
S_\phi^{X^\epsilon}(t)<0
\right\}
\subseteq
\left\{
(X^\epsilon,t)\not\models\phi
\right\}
\subseteq
\left\{
S_\phi^{X^\epsilon}(t)\leq0
\right\}.
\label{eq:violation_sandwich_LDP}
\end{align}
Applying the lower and upper bounds of the LDP in
Corollary~\ref{cor:fixed_time_robustness_LDP} to the two outer events
in \eqref{eq:violation_sandwich_LDP} yields
\eqref{eq:rare_violation_bounds}. If
\eqref{eq:boundary_rate_condition} holds, then the lower and upper
large-deviation bounds coincide, which gives
\eqref{eq:exact_failure_exponent}.
\end{proof}

Condition \eqref{eq:boundary_rate_condition} states that approaching the
logical boundary does not offer an exponentially cheaper mechanism than
entering the strict violation region. When this condition holds,
$I_{\phi,t}^{\mathrm{fail}}$ has a direct interpretation: it is the
minimum pathwise cost required to produce a temporal-logic failure. A completely analogous statement can be formulated for rare satisfaction events using the sets $\{\rho_\phi>0\}$ and $\{\rho_\phi\geq0\}$.

\subsection{Small-noise diffusions}
The preceding results are particularly natural for small-noise
diffusion processes. Consider
\begin{align}
dX_t^\epsilon
=
b(X_t^\epsilon)\,dt
+
\sqrt{\epsilon}\,
\sigma(X_t^\epsilon)\,dW_t,
\qquad
X_0^\epsilon=x_0,
\label{eq:FW_diffusion}
\end{align}
on the interval $[0,T+H_\phi]$, where $W$ is a standard Brownian motion in $\mathbb{R}^m$, $b : \mathbb{R}^d \to \mathbb{R}^d$, and $\sigma : \mathbb{R}^d\to\mathbb{R}^{d\times m}$. Under standard Freidlin--Wentzell regularity assumptions, the family $\{X^\epsilon\}_{\epsilon>0}$ satisfies a path-space LDP as
$\epsilon\downarrow0$; see
\citep{FreidlinWentzell2012}. For simplicity, suppose that $b$ and
$\sigma$ are bounded and Lipschitz and that
\begin{align}
a(x)
:=
\sigma(x)\sigma(x)^\top
\end{align}
is uniformly positive definite. Then the corresponding rate function is
\begin{align}
I_{T+H_\phi}(x)
=
\begin{cases}
\displaystyle
\frac12
\int_0^{T+H_\phi}
\big(
\dot{x}(s)-b(x(s))
\big)^\top
a(x(s))^{-1}
\big(
\dot{x}(s)-b(x(s))
\big)
\,ds,
\\
\hfill
x\text{ absolutely continuous and }x(0)=x_0,
\\[3mm]
+\infty,
\qquad\text{otherwise}.
\end{cases}
\label{eq:FW_rate_function}
\end{align}
Combining the Freidlin--Wentzell LDP \citep{dos2019freidlin} with
Theorem~\ref{thm:TLOP_LDP} gives the following consequence.

\begin{corollary}
\label{cor:FW_TLOP}
Consider the small-noise diffusion \eqref{eq:FW_diffusion} and let $\phi$ satisfy the assumptions of
Proposition~\ref{prop:robustness_map_continuity}. Then
$\{S_\phi^{X^\epsilon}\}_{\epsilon>0}$ satisfies an LDP on
$C([0,T],\mathbb{R})$ with rate function
\begin{align}
J_{\phi,T}(z)
=
\inf_{\mathcal{R}_{\phi,T}(x)=z}
I_{T+H_\phi}(x),
\label{eq:FW_TLOP_rate}
\end{align}
\end{corollary}
\noindent The action functional \eqref{eq:FW_rate_function} has a useful control interpretation. The deterministic zero-noise system satisfies
\begin{align}
\dot{x}(t)=b(x(t)).
\end{align}
A rare path $x$ deviates from this nominal dynamics. Formally, such a
deviation may be generated by introducing an auxiliary input $u$ through
\begin{align}
\dot{x}(t)
=
b(x(t))
+
\sigma(x(t))u(t).
\label{eq:controlled_FW}
\end{align}
The corresponding large-deviation action represents the minimum
quadratic energy of such an input. Hence the most likely rare temporal
failure can be interpreted as a minimum-energy trajectory capable of
driving the deterministic system into the temporal-logic violation set. In particular, whenever
\eqref{eq:boundary_rate_condition} holds,
\begin{align}
\label{eq:temporal_failure_action}
I_{\phi,t}^{\mathrm{fail}}
=
\inf_{\rho_\phi(x,t)\leq0}
I_{T+H_\phi}(x).
\end{align}
The quantity $I_{\varphi,t}^{\mathrm{fail}}$
characterizes the minimum large-deviation action required to produce a temporal-logic violation. Under \eqref{eq:controlled_FW}, and assuming
$I_{\varphi,t}^{\mathrm{fail}}<\infty$,
the probability of temporal-logic violation satisfies the logarithmic asymptotic relation
\begin{equation}
\label{eq:temporal_failure_probability}
\lim_{\varepsilon\downarrow 0}
\varepsilon\log
\mathbb{P}
\left(
(X^\varepsilon,t)\not\models\varphi
\right)
=
-I_{\varphi,t}^{\mathrm{fail}}
\iff
\mathbb{P}
\left(
(X^\varepsilon,t)\not\models\varphi
\right)
=
\exp\left(
-\frac{
I_{\varphi,t}^{\mathrm{fail}}+o(1)
}{\varepsilon}
\right),
\qquad \varepsilon\downarrow 0.
\end{equation}
In particular, when
$I_{\varphi,t}^{\mathrm{fail}}>0$,
temporal-logic violation becomes exponentially unlikely as the noise intensity $\varepsilon$ tends to zero. The large-deviation rate determines the leading exponential decay, while sub-exponential factors are not characterized by this result. If a minimizer in \eqref{eq:temporal_failure_action} exists, then 
\begin{align}
x_\phi^\star
\in
\arg\min_{\rho_\phi(x,t)\leq0}
I_{T+H_\phi}(x),
\label{eq:most_likely_failure_path}
\end{align}
where $x_\phi^\star$ may be interpreted as a most likely trajectory by
which the small-noise stochastic system approaches temporal failure at
the exponential scale.
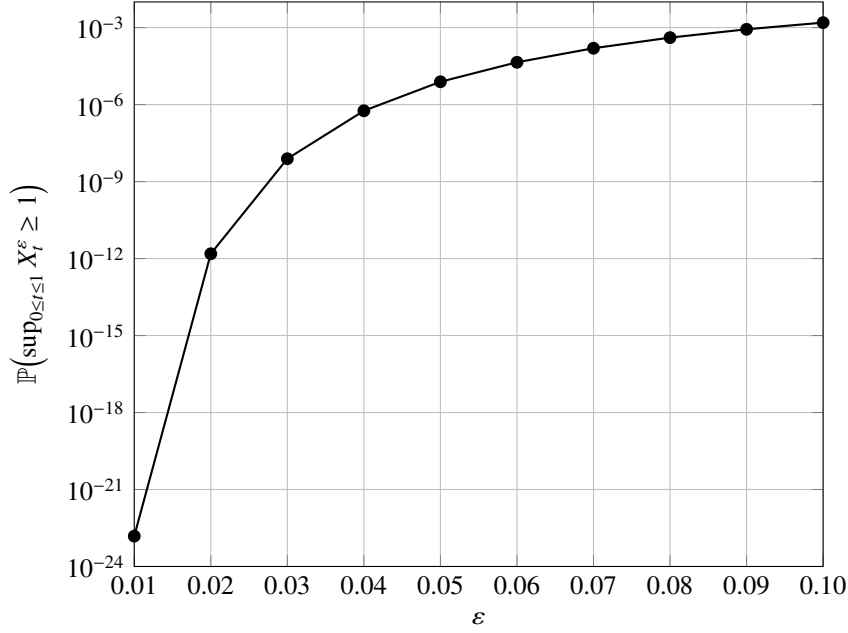
\begin{figure}[t]
    \centering
    \begin{tikzpicture}
        \begin{semilogyaxis}[
            width=0.65\linewidth,
            height=0.55\linewidth,
            xlabel={$\varepsilon$},
            ylabel={
                $\mathbb{P}\!\left(
                \sup_{0\leq t\leq 1} X_t^\varepsilon \geq 1
                \right)$
            },
            xmin=0.01,
            xmax=0.10,
            xtick={0.01,0.02,0.03,0.04,0.05,0.06,0.07,0.08,0.09,0.10},
            xticklabels={$0.01$,$0.02$,$0.03$,$0.04$,$0.05$, $0.06$,$0.07$,$0.08$,$0.09$,$0.10$},
            ymin=1e-24,
            ymax=1e-2,
            grid=both,
            mark size=2pt,
        ]

        \addplot[
            thick,
            mark=*,
        ]
        coordinates {
            (0.01,1.5239706e-23)
            (0.02,1.5374598e-12)
            (0.03,7.7640365e-09)
            (0.04,5.7330314e-07)
            (0.05,7.7442164e-06)
            (0.06,4.4557091e-05)
            (0.07,1.5705228e-04)
            (0.08,4.0695202e-04)
            (0.09,8.5812067e-04)
            (0.10,1.5654023e-03)
        };

        \end{semilogyaxis}
    \end{tikzpicture}
    \caption{
        Probability of reaching the level $K=1$ within the time horizon
        $T=1$ for the small-noise Brownian process
        $X_t^\varepsilon=\sqrt{\varepsilon}W_t$.
        The exact probability is
        $2\left[1-\Phi(1/\sqrt{\varepsilon})\right]$.
        The logarithmic vertical axis highlights the rapid decay of the
        target-attainment probability as $\varepsilon\downarrow0$.
    }
    \label{fig:small_noise_reaching_probability}
\end{figure}
\subsection{Consistency with the Brownian target-attainment example}
Returning to the target-attainment example in \eqref{eq:small_noise_BM}, we have $b(x)=0, \sigma(x)=1$. The Freidlin--Wentzell rate function therefore reduces to
\[
I_T(x)
=
\frac12\int_0^T|\dot{x}(s)|^2\,ds.
\]
For $\phi = \lozenge_{[0,T]}(X\geq K)$, the minimum action required for satisfaction is
\[
\inf_{\sup_{0\leq s\leq T}x(s)\geq K}
I_T(x)
=
\frac{(K-x_0)^2}{2T},
\]
with minimizing path given by $x^\star(s) = x_0 + \frac{K-x_0}{T}s$. The resulting large-deviation prediction is therefore
\begin{align*}
\lim_{\epsilon\downarrow0}
\epsilon
\log
\mathbb{P}\big((X^\epsilon,0)\models\phi\big)
=
-\frac{(K-x_0)^2}{2T},
\end{align*}
which agrees exactly with the exponent obtained from the explicit
reflection-principle probability
\eqref{eq:exact_small_noise_probability}. This example illustrates the advantage of the large-deviation
formulation. For Brownian motion the temporal-event probability can be
computed explicitly. For more general stochastic processes such closed
forms are rarely available. The path-space formulation instead reduces
the logarithmic rare-event asymptotics to a variational problem over
trajectories, while Theorem~\ref{thm:TLOP_LDP} transfers this structure
directly to the associated temporal-logical observable process.

\subsection{A Two-Stage Temporal Rare Event:
Target Attainment Followed by a Downward Excursion}

The preceding example considered a single target-attainment
specification. We now investigate a more involved temporal
event containing two distinct temporal requirements. This
example illustrates how the structure of a temporal formula
affects both the large-deviation rate and the minimum-action
trajectory associated with a rare event. Consider the small-noise Brownian process as in \eqref{eq:small_noise_BM} given by
\begin{equation}
\label{eqn_two_stage_rare_event_process}
X_t^\varepsilon=x_0+\sqrt{\varepsilon}W_t,
\qquad t\geq 0.
\end{equation}
The process may be interpreted as a stylized model of fluctuations in the logarithm of an asset price. Let $L<x_0<K$ and $0<T_1<T_2$. Suppose that an investor is interested in the following temporal event: the process
must first reach an upper target $K$ before time $T_1$ and subsequently reach a lower level $L$ during the interval $[T_1,T_2]$. This requirement is expressed through the temporal formula
\begin{equation}
\label{eqn_two_stage_rare_event_formula}
\varphi :=
\Diamond_{[0,T_1]}(X\geq K)
\wedge
\Diamond_{[T_1,T_2]}(X\leq L).
\end{equation}
The atomic predicate functions corresponding to the upper and lower thresholds are $h_{\mu_K}(x) = x-K$, and $h_{\mu_L}(x) = L-x$ respectively. Both predicate functions are globally Lipschitz with Lipschitz constant one. The temporal-logical observable at time $t=0$ is 
\begin{equation}
\begin{aligned}
S_\varphi^{X^\varepsilon}(0)
=
\min\Big\{
\sup_{s\in[0,T_1]}(X_s^\varepsilon-K),
\sup_{s\in[T_1,T_2]}(L-X_s^\varepsilon)
\Big\}.
\end{aligned}
\end{equation}
Since the sample paths are continuous and the formula contains only threshold predicates, conjunction, and eventually operators, Boolean satisfaction is equivalent to non-negative robustness in this example. Therefore,
\begin{equation}
\begin{aligned}
p_\varepsilon
:=
\mathbb P((X^\varepsilon,0)\models\varphi)
=
\mathbb P\left(
S_\varphi^{X^\varepsilon}(0)\geq0
\right).
\end{aligned}
\end{equation}
We seek to characterize the exponential decay of $p_\varepsilon$ as $\varepsilon \downarrow0$. By Schilder's theorem, the family
$\{X^\varepsilon\}_{\varepsilon>0}$ satisfies an LDP on $C([0,T_2],\mathbb R)$ with speed $\varepsilon^{-1}$ and good rate function
\begin{equation}
I_{T_2}(x)=
\begin{cases}
\displaystyle
\frac12\int_0^{T_2}|\dot x(s)|^2\,ds,
&
x\text{ absolutely continuous}, 
x(0)=x_0, \\
+\infty,&\text{otherwise}.
\end{cases}
\end{equation}
Define the set of trajectories satisfying the temporal specification as
\begin{equation}
\begin{aligned}
\mathcal A_\varphi:=
\Big\{x\in C([0,T_2],\mathbb R) \mid
\sup_{s\in[0,T_1]}x(s)\geq K,
\inf_{s\in[T_1,T_2]}x(s)\leq L
\Big\}.
\end{aligned}
\end{equation}
The minimum action associated with this event is $I_\varphi^\star := \inf_{x\in\mathcal A_\varphi}I_{T_2}(x)$. The following proposition provides an explicit
characterization of this minimum action and
the corresponding logarithmic rare-event probability.

\begin{proposition}
Given the small-noise process \eqref{eqn_two_stage_rare_event_process} and the formula $\varphi$ defined using \eqref{eqn_two_stage_rare_event_formula}, let $a:=K-x_0>0$, and $b:=K-L>0$. Define $\tau_0:=\frac{aT_2}{a+b}$, and $\tau^\star:=\min\{T_1,\tau_0\}$. Then, 
\begin{align}
\lim_{\varepsilon\downarrow0}
\varepsilon\log
p_\varepsilon
&=
-I_\varphi^\star, \quad \text{where}, \\
I_\varphi^\star
=
\frac{a^2}{2\tau^\star}
+
\frac{b^2}{2(T_2-\tau^\star)} 
&=
\begin{cases}
\displaystyle
\frac{(a+b)^2}{2T_2},
&T_1\geq\tau_0,\\
\displaystyle
\frac{a^2}{2T_1}
+
\frac{b^2}{2(T_2-T_1)},
&T_1<\tau_0.
\end{cases}
\end{align}
Further, the minimum-action trajectory is given by
\begin{equation}
x^\star(s)=
\begin{cases}
\displaystyle
x_0+\frac{K-x_0}{\tau^\star}s,
&0\leq s\leq\tau^\star,\\[12pt]
\displaystyle
K-\frac{K-L}{T_2-\tau^\star}
(s-\tau^\star),
&\tau^\star<s\leq T_2.
\end{cases}
\end{equation}
\end{proposition}

\begin{proof}
We first determine the minimum action among the trajectories satisfying $\varphi$. Let $x\in\mathcal A_\varphi$ be an absolutely continuous trajectory starting at $x_0$. By the definition of $\mathcal A_\varphi$, there exist times $\tau\in(0,T_1], v\in[T_1,T_2]$, such that $x(\tau)\geq K, x(v)\leq L$. Since $K>L$ and $\tau\leq T_1\leq v$,
necessarily $\tau<v$. Using the familiar machinery by invoking the fundamental theorem of calculus, followed by applying the Cauchy--Schwarz inequality and using the fact that $x(\tau)\geq K>x_0$, we obtain
\begin{align*}
x(\tau)-x_0
= 
\int_0^\tau\dot x(s)\,ds 
\implies 
|x(\tau)-x_0|^2
\leq
\tau\int_0^\tau|\dot x(s)|^2\,ds 
\implies 
\frac12\int_0^\tau|\dot x(s)|^2\,ds
\geq
\frac{(K-x_0)^2}{2\tau}
=
\frac{a^2}{2\tau}.
\end{align*}
Similarly, use the fact that $x(\tau)\geq K$ and $x(v)\leq L$,
we have $|x(v)-x(\tau)|\geq K-L=b$. Then, using similar machinery, we get
\begin{align*}
x(v)-x(\tau)
=
\int_\tau^v\dot x(s)\,ds 
\implies 
|x(v)-x(\tau)|^2
\leq
(v-\tau)
\int_\tau^v|\dot x(s)|^2\,ds 
\implies 
\frac12\int_\tau^v|\dot x(s)|^2\,ds
\geq
\frac{b^2}{2(v-\tau)}.
\end{align*}
Combining these inequalities and using
$v\leq T_2$ yields
\begin{align*}
I_{T_2}(x)
=
\frac12\int_0^{T_2}|\dot x(s)|^2\,ds
\geq
\frac{a^2}{2\tau}
+
\frac{b^2}{2(v-\tau)}
\geq
\frac{a^2}{2\tau}
+
\frac{b^2}{2(T_2-\tau)}.
\end{align*}
Hence every admissible trajectory satisfies
\begin{equation*}
I_{T_2}(x)
\geq
\inf_{0<\tau\leq T_1}
\left\{
\frac{a^2}{2\tau}
+
\frac{b^2}{2(T_2-\tau)}
\right\}.
\end{equation*}
For $0<\tau\leq T_1 < T_2$, define
\begin{align*}
f(\tau)
:=
\frac{a^2}{2\tau}
+
\frac{b^2}{2(T_2-\tau)}
\implies
f'(\tau)
=
-\frac{a^2}{2\tau^2}
+
\frac{b^2}{2(T_2-\tau)^2} 
\implies
f''(\tau)
=
\frac{a^2}{\tau^3}
+
\frac{b^2}{(T_2-\tau)^3}>0.
\end{align*}
Thus $f$ is strictly convex. Then, $f'(\tau) = 0 \implies \frac{a}{\tau} = \frac{b}{T_2-\tau}$. Solving for $\tau$ gives $\tau_0=\frac{aT_2}{a+b}$. Since the first temporal requirement imposes $\tau\leq T_1$, the constrained minimizer is $\tau^\star=\min\{T_1,\tau_0\}$. So, 
\begin{equation*}
I_{T_2}(x)
\geq
\frac{a^2}{2\tau^\star}
+
\frac{b^2}{2(T_2-\tau^\star)}.
\end{equation*}
To show that this lower bound is attainable, consider the piecewise linear trajectory
\begin{equation*}
x^\star(s)=
\begin{cases}
\displaystyle
x_0+\frac{a}{\tau^\star}s,
&0\leq s\leq\tau^\star,\\
\displaystyle
K-\frac{b}{T_2-\tau^\star}
(s-\tau^\star),
&\tau^\star<s\leq T_2.
\end{cases}
\end{equation*}
By construction, $x^\star(0)=x_0, x^\star(\tau^\star)=K, x^\star(T_2)=L$. Since $\tau^\star\leq T_1<T_2$, this trajectory satisfies both temporal
requirements and belongs to
$\mathcal A_\varphi$. Its action is
\begin{align*}
I_{T_2}(x^\star)
=
\frac12\int_0^{\tau^\star}
\left(\frac{a}{\tau^\star}\right)^2ds
+
\frac12\int_{\tau^\star}^{T_2}
\left(\frac{b}{T_2-\tau^\star}\right)^2ds
=
\frac{a^2}{2\tau^\star}
+
\frac{b^2}{2(T_2-\tau^\star)}.
\end{align*}
Thus the lower bound is attained, proving the claimed formula for $I_\varphi^\star$. If $T_1\geq\tau_0$, substitution of $\tau^\star=\tau_0$ gives $I_\varphi^\star = \frac{(a+b)^2}{2T_2}$. If $T_1<\tau_0$, then
$\tau^\star=T_1$, giving $I_\varphi^\star = \frac{a^2}{2T_1} + \frac{b^2}{2(T_2-T_1)}$. It remains to establish the logarithmic rare-event probability. Since the maximum and minimum functionals are continuous under the uniform norm, $\mathcal A_\varphi$ is closed in $C([0,T_2],\mathbb R)$. Schilder's large-deviation upper bound gives
\begin{equation}
\label{eqn_Schilder_upper_bound}
\limsup_{\varepsilon\downarrow0}
\varepsilon\log p_\varepsilon
\leq-I_\varphi^\star.
\end{equation}
For the lower bound, consider the open set
\begin{equation*}
\begin{aligned}
\mathcal A_\varphi^\circ
:=
\Big\{x \mid
\sup_{s\in[0,T_1]}x(s)>K,
\inf_{s\in[T_1,T_2]}x(s)<L
\Big\}.
\end{aligned}
\end{equation*}
For every $\delta>0$, construct a piecewise linear trajectory $x_\delta$ joining $x_0$ at time $0$ to $K+\delta$
at time $\tau^\star$, and subsequently
to $L-\delta$ at time $T_2$. Then $x_\delta\in\mathcal A_\varphi^\circ$, and $I_{T_2}(x_\delta) \longrightarrow
I_\varphi^\star$ as $\delta \downarrow 0$. Hence,
\begin{equation*}
\inf_{x\in\mathcal A_\varphi^\circ}
I_{T_2}(x)
=
\inf_{x\in\mathcal A_\varphi}
I_{T_2}(x)
=
I_\varphi^\star.
\end{equation*}
Since $\mathcal A_\varphi^\circ
\subseteq\mathcal A_\varphi$, the large-deviation lower bound yields
\begin{equation}
\label{eqn_Schilder_lower_bound}
\liminf_{\varepsilon\downarrow0}
\varepsilon\log p_\varepsilon
\geq-I_\varphi^\star.
\end{equation}
Combining the upper bound \eqref{eqn_Schilder_upper_bound} and the lower bound \eqref{eqn_Schilder_lower_bound}, we get the equality given by $\lim_{\varepsilon\downarrow0}
\varepsilon\log p_\varepsilon = -I_\varphi^\star$. 
\end{proof}
The preceding proposition shows that the
large-deviation rate is explicitly affected by the temporal deadline $T_1$.
When $T_1\geq\tau_0$, the first deadline
does not constrain the minimum-action trajectory, and the rare-event rate is $I_\varphi^\star = \frac{(a+b)^2}{2T_2}$. In this regime, the upward and downward portions of the minimum-action trajectory have equal absolute slopes. On the other hand, when $T_1<\tau_0$, the first deadline becomes active. The trajectory must reach $K$ precisely at $T_1$, requiring a larger initial rate of change and an increased minimum action. For example, let $x_0=0, K=1, L=-1, T_2=2$. Then $a=1$, $b=2$, and $\tau_0=2/3$. For $T_1=1$, the first deadline is inactive, giving $I_\varphi^\star=\frac94$. If the deadline is tightened to $T_1=1/2$, it becomes active, and $I_\varphi^\star = \frac{1}{2(1/2)} + \frac{4}{2(2-1/2)} = \frac73$. Thus, tightening the temporal requirement increases the large-deviation rate, making the specified event exponentially less likely in the small-noise limit. This example demonstrates how a compound temporal-logic specification induces a constrained minimum-action problem whose solution depends explicitly on the temporal parameters. In contrast with the single-barrier example, the logical structure now determines both the optimal intermediate hitting time and the geometry of the most likely rare-event trajectory. The same variational perspective extends to more general temporal specifications whenever a suitable path-space LDP and continuity of the robustness map are available.

\section{Conclusion and Future Outlook} \label{sec_11}
In this manuscript, we developed a measure-theoretic framework for TLSP, wherein temporal logic specifications induce stochastic observable processes that quantify temporal satisfaction and robustness. We established the measurability, algebraic structure, and sample-path regularity of these observables, and characterized
their connections to extremal functionals and first hitting times of stochastic processes. Furthermore, we derived pathwise perturbation bounds and identified the role of probability mass near the logical boundary in determining the stability of temporal satisfaction
probabilities. By introducing a topology on the observable-process space and a semantic pseudometric on temporal formulas, we established convergence results for TLOPs under perturbations of temporal specifications, together with corresponding
distributional convergence results. We also established a process-level large-deviation principle for TLOPs, providing a variational characterization of rare
temporal events through minimum-action trajectories.
The Brownian examples, including the two-stage temporal
rare-event problem, illustrate how temporal constraints
influence the associated large-deviation rates and
minimum-action trajectories.

Several research directions arise naturally from this
framework. First, the connection between temporal satisfaction probabilities and Doob $h$-transforms could be investigated by constructing suitable space-time harmonic functions from conditional temporal satisfaction probabilities, potentially using augmented states to account for path-dependent specifications, thereby enabling stochastic dynamics conditioned on temporal events. Second, Feynman--Kac representations could provide a connection between temporal-logic observables and partial differential equations through suitable path-dependent costs or softened temporal penalties, offering analytical approaches to evaluating
temporally structured expectations. Third, conditioned diffusions could be studied to characterize the stochastic trajectories of a system conditioned on satisfying a given temporal specification, and to relate their conditional path distributions to the minimum-action trajectories obtained through large-deviation analysis. Fourth, Schr\"odinger bridges could provide an
entropy-minimizing framework for modifying stochastic
path distributions subject to temporal-logic constraints,
potentially yielding probabilistic mechanisms for
generating trajectories that satisfy prescribed
temporal requirements. Finally, extending the framework to richer temporal operators, particularly the \emph{until} operator, and to stochastic processes with discontinuous sample paths would require further investigation of nested temporal extrema, measurability, and continuity under appropriate path-space topologies. These directions suggest that TLSP may provide a
common mathematical setting for investigating temporal
logic through stochastic analysis, connecting logical
specifications with probability laws, stochastic
conditioning, and rare-event behaviour.


\section*{Declaration of generative AI and AI-assisted technologies in the manuscript preparation process}
During the preparation of this manuscript, the author used OpenAI’s ChatGPT to polish and improve its clarity of the content. The author reviewed the output and edited it needed and shall take the full responsibility for the content of the published article.

\bibliographystyle{cas-model2-names}
\bibliography{references}

@book{oksendal2003stochastic,
  title={Stochastic Differential Equations},
  author={{\O}ksendal, B.},
  year={2003},
  publisher={Springer}
}

@article{sleiman_acc,
  author={Safaoui, Sleiman and Lindemann, Lars and Shames, Iman and Summers, Tyler H.},
  booktitle={2022 American Control Conference (ACC)}, 
  title={Risk-Bounded Temporal Logic Control of Continuous-Time Stochastic Systems}, 
  year={2022},
  volume={},
  number={},
  pages={1555-1562},
  doi={10.23919/ACC53348.2022.9867734}}

@book{karatzas1991brownian,
  title={Brownian Motion and Stochastic Calculus},
  author={Karatzas, I. and Shreve, S.},
  year={1991},
  publisher={Springer}
}

@article{pnueli1977temporal,
  title={The Temporal Logic of Programs},
  author={Pnueli, A.},
  journal={FOCS},
  year={1977}
}

@book{baier2008principles,
  title={Principles of Model Checking},
  author={Baier, C. and Katoen, J.-P.},
  year={2008},
  publisher={MIT Press}
}

@article{maler2004monitoring,
  title={Monitoring Temporal Properties of Continuous Signals},
  author={Maler, O. and Nickovic, D.},
  journal={FORMATS},
  year={2004}
}

@article{donze2010robust,
  title={Robust Satisfaction of Temporal Logic over Real-Valued Signals},
  author={Donz{\'e}, A. and Maler, O.},
  journal={HSCC},
  year={2010}
}

@article{kwiatkowska2007stochastic,
  title={Stochastic Model Checking},
  author={Kwiatkowska, M. et al.},
  journal={Formal Methods},
  year={2007}
}

@article{abate2008probabilistic,
  title={Probabilistic reachability and safety for controlled discrete time stochastic hybrid systems},
  author={Abate, Alessandro and Prandini, Maria and Lygeros, John and Sastry, Shankar},
  journal={Automatica},
  volume={44},
  number={11},
  pages={2724--2734},
  year={2008},
  publisher={Elsevier}
}

@book{billingsley1995probability,
  title={Probability and Measure},
  author={Billingsley, P.},
  year={1995},
  publisher={Wiley}
}

@book{durrett2019probability,
  title={Probability: Theory and Examples},
  author={Durrett, R.},
  year={2019},
  publisher={Cambridge}
}

@book{shiryaev1996probability,
  title={Probability},
  author={Shiryaev, A.},
  year={1996},
  publisher={Springer}
}

@book{ethier2009markov,
  title={Markov Processes},
  author={Ethier, S. and Kurtz, T.},
  year={2009},
  publisher={Wiley}
}

@book{revuz1999continuous,
  title={Continuous Martingales and Brownian Motion},
  author={Revuz, D. and Yor, M.},
  year={1999},
  publisher={Springer}
}

@article{uhlenbeck1930theory,
  title={On the theory of the Brownian motion},
  author={Uhlenbeck, George E and Ornstein, Leonard S},
  journal={Physical review},
  volume={36},
  number={5},
  pages={823},
  year={1930},
  publisher={APS}
}

@book{vershynin2018high,
  title={High-dimensional probability: An introduction with applications in data science},
  author={Vershynin, Roman},
  volume={47},
  year={2018},
  publisher={Cambridge university press}
}

@article{Madsen2018Metrics,
  author    = {Curtis Madsen and
               Prashant Vaidyanathan and
               Sadra Sadraddini and
               Cristian-Ioan Vasile and
               Nicholas A. DeLateur and
               Ron Weiss and
               Douglas Densmore and
               Calin Belta},
  title     = {Metrics for Signal Temporal Logic Formulae},
  booktitle = {2018 IEEE Conference on Decision and Control (CDC)},
  pages     = {1542--1547},
  year      = {2018},
  doi       = {10.1109/CDC.2018.8619541}
}

@article{fainekos2009robustness,
  title={Robustness of temporal logic specifications for continuous-time signals},
  author={Fainekos, Georgios E and Pappas, George J},
  journal={Theoretical Computer Science},
  volume={410},
  number={42},
  pages={4262--4291},
  year={2009},
  publisher={Elsevier}
}

@article{fasoni2020robustness,
  title={Robustness of temporal logic specifications},
  author={Fainekos, Georgios E and Pappas, George J},
  journal={International Workshop on Formal Approaches to Software Testing},
  pages={178--192},
  year={2006},
  publisher={Springer}
}

@article{Bartocci2015,
  author  = {Ezio Bartocci and Luca Bortolussi and Laura Nenzi and Guido Sanguinetti},
  title   = {System Design of Stochastic Models Using Robustness of Temporal Properties},
  journal = {Theoretical Computer Science},
  volume  = {587},
  pages   = {3--25},
  year    = {2015}
}

@inproceedings{Lindemann2021,
  author    = {Lars Lindemann and Nikolai Matni and George J. Pappas},
  title     = {{STL} Robustness Risk over Discrete-Time Stochastic Processes},
  booktitle = {Proceedings of the 60th IEEE Conference on Decision and Control (CDC)},
  pages     = {1329--1335},
  year      = {2021},
  publisher = {IEEE}
}

@article{lindemann2023risk,
  title={Risk of stochastic systems for temporal logic specifications},
  author={Lindemann, Lars and Jiang, Lejun and Matni, Nikolai and Pappas, George J},
  journal={ACM Transactions on Embedded Computing Systems},
  volume={22},
  number={3},
  pages={1--31},
  year={2023},
  publisher={ACM New York, NY}
}

@article{Frankova2019,
  author  = {Fra\v{n}kov\'a, Dana},
  title   = {Regulated Functions with Values in Banach Space},
  journal = {Mathematica Bohemica},
  volume  = {144},
  number  = {4},
  pages   = {437--456},
  year    = {2019},
  doi     = {10.21136/MB.2019.0124-19}
}

@article{dos2019freidlin,
  title={Freidlin--Wentzell {LDP} in path space for McKean--Vlasov equations and the functional iterated logarithm law},
  author={Dos Reis, Gon{\c{c}}alo and Salkeld, William and Tugaut, Julian},
  journal={The Annals of Applied Probability},
  year={2019},
  volume={29},
  number={3},
  pages={1487--1540},
}

@article{cleaveland2022risk,
  title={Risk verification of stochastic systems with neural network controllers},
  author={Cleaveland, Matthew and Lindemann, Lars and Ivanov, Radoslav and Pappas, George J},
  journal={Artificial Intelligence},
  volume={313},
  pages={103782},
  year={2022},
  publisher={Elsevier}
}

@book{DemboZeitouni2010,
  author    = {Amir Dembo and Ofer Zeitouni},
  title     = {Large Deviations Techniques and Applications},
  series    = {Stochastic Modelling and Applied Probability},
  volume    = {38},
  edition   = {2},
  publisher = {Springer},
  year      = {2010},
  doi       = {10.1007/978-3-642-03311-7}
}

@book{FreidlinWentzell2012,
  author    = {Mark I. Freidlin and Alexander D. Wentzell},
  title     = {Random Perturbations of Dynamical Systems},
  series    = {Grundlehren der mathematischen Wissenschaften},
  volume    = {260},
  edition   = {3},
  publisher = {Springer},
  year      = {2012},
  doi       = {10.1007/978-3-642-25847-3}
}

\end{document}